\documentclass{amsart}
\usepackage{amssymb,amsmath, amsthm,latexsym}
\usepackage{graphics}
\usepackage{amscd}
\newcommand{\cal}[1]{\mathcal{#1}}
\theoremstyle{plain}
\newtheorem{theorem}{Theorem}
\newtheorem{coro}{Corollary}
\newtheorem{lemma}{Lemma}[section]
\newtheorem{theo}[lemma]{Theorem}
\newtheorem{proposition}[lemma]{Proposition}
\newtheorem{corollary}[lemma]{Corollary}
\theoremstyle{definition}
\newtheorem{definition}[lemma]{Definition}
\parskip=\bigskipamount

\let\egthree=\phi
\let\phi=\varphi
\let\varphi=\egthree

\begin{document}
\title{Isometry groups of proper ${\rm CAT}(0)$-spaces}
\author{Ursula Hamenst\"adt}

\date{February 10, 2008}

\begin{abstract}
Let $X$ be a proper ${\rm CAT}(0)$-space
and let $G$
be a closed subgroup of the isometry group ${\rm Iso}(X)$ of $X$.
We show that if $G$ is non-elementary and 
contains a rank-one element then its second
bounded cohomology group with coefficients in the regular
representation is non-trivial. As a consequence,
up to passing to an open subgroup of finite index,
either $G$ is a compact extension of a totally
disconnected group or $G$ is a compact extension
of a simple Lie group of rank one.
\end{abstract}

\maketitle

\section{Introduction}

A geodesic metric space $(X,d)$ is called \emph{proper} if
closed balls in $X$ of finite radius are compact. 
A proper ${\rm CAT}(0)$-metric space $X$
can be compactified by adding the
\emph{visual boundary} $\partial X$.
The
isometry group ${\rm Iso}(X)$ of $X$, equipped with the compact
open topology, is a locally compact $\sigma$-compact
topological group which acts as a group of homeomorphisms
on $\partial X$. The \emph{limit set} $\Lambda$ of a subgroup
$G$ of ${\rm Iso}(X)$ is the set of accumulation points
in $\partial X$ of an orbit of the action of $G$ on $X$. 
The group $G$ is called
\emph{elementary} if either its limit set consists of at most
two points or if 
$G$ fixes a point in $\partial X$.

For every $g\in {\rm Iso}(X)$ the \emph{displacement function}
of $g$ is the function $x\to d(x,gx)$. The isometry $g$ is called
\emph{semisimple} if its displacement function assumes a
minimum on $X$. If this minimum vanishes
then $g$ has a fixed point in $X$ and is called \emph{elliptic},
and otherwise $g$ is called \emph{axial}.
If $g$ is axial then the closed
convex subset $A$ of $X$ on which the displacement
function is minimal is isometric to a product 
space $C\times\mathbb{R}$ where
$g$ acts on 
each of the geodesics $\{x\}\times \mathbb{R}$ as a translation.
Such a geodesic is called an \emph{axis} for $g$. We refer
to the books \cite{B95,BGS85,BH99} for basic properties of 
${\rm CAT}(0)$-spaces and for references.
 
Call an axial isometry $g$ of $X$ \emph{rank-one} if there
is an axis $\gamma$ for $g$ which does not bound a flat half-plane.
Here by a flat half-plane we mean 
a totally geodesic embedded isometric copy of an euclidean half-plane
in $X$.

A \emph{compact extension} of a topological group $H$
is a topological group $G$ which contains a compact
normal subgroup $K$ such that $H=G/K$ as topological groups.
Extending earlier results for isometry groups of proper hyperbolic
geodesic metric spaces \cite{H08b,MMS04,MS04} we show.

\begin{theorem}\label{thm1} 
Let $X$ be a proper ${\rm CAT}(0)$-space
and let $G< {\rm Iso}(X)$ be a
closed subgroup. Assume that $G$ is non-elementary
and contains a rank-one element. 
Then one of the following two
possibilities holds.
\begin{enumerate}
\item Up to passing to an open subgroup of finite
index, $G$ is a compact extension of a simple Lie group of rank one.
\item $G$ is a compact extension of a totally disconnected group.
\end{enumerate}
\end{theorem}

Caprace and Monod 
showed the following version of Theorem 1 (Corollary 1.7 of
\cite{CM08}): A ${\rm CAT}(0)$-space $X$ is called 
\emph{irreducible} if it is not a non-trivial metric product.
Let $X\not=\mathbb{R}$ be an irreducible proper
${\rm CAT}(0)$-space with finite dimensional Tits boundary.
Assume that the isometry group ${\rm Iso}(X)$ of $X$
does not have a global fixed point in $\partial X$ and that
its action on $X$ does not preserve a non-trivial
closed convex subset of $X$. 
Then ${\rm Iso}(X)$ is either totally disconnected
or an almost connected simple Lie group with trivial center.

We also note the following consequence (see Corollary 1.24 of
\cite{CM08}). 

\begin{coro}\label{rank}
Let $M$ be a closed Riemannian manifold of non-positive
sectional curvature. If the universal
covering $\tilde M$ of $M$ is irreducible 
and if the isometry group of 
$\tilde M$ contains a 
parabolic element then $M$ is locally symmetric.
\end{coro}

Our proof of Theorem 1 is different
from the approach of Caprace and Monod and 
uses second bounded
cohomology for locally compact topological groups $G$
with coefficients in a 
\emph{Banach module} for $G$. Such a Banach module 
is a separable Banach space $E$ together with a
continuous
homomorphism of $G$ into the group of linear
isometries of $E$. For every such Banach module $E$
for $G$ and every $i\geq 1$, the
group $G$ naturally acts on the vector space
$C_b(G^i,E)$ of continuous bounded maps $G^i\to E$.
If we denote by
$C_b(G^i,E)^G\subset C_b(G^i,E)$
the linear subspace of all $G$-invariant such maps, then the
\emph{second continuous bounded
cohomology group} $H_{cb}^2(G,E)$ of $G$
with coefficients $E$
is defined as the second cohomology group of the complex
\begin{equation} 0\to
C_b(G,E)^G \xrightarrow{d} C_b
(G^2,E)^G\xrightarrow{d} \dots \notag
\end{equation} with the usual
homogeneous coboundary operator $d$ (see \cite{M}).
We write $H_b^2(G,E)$ to denote the second continuous
bounded cohomology group of $G$ with coefficients
$E$ where $G$ is equipped with the discrete
topology (but where the action of $G$ on $E$ is defined by
a continuous homomorphism of $G$ equipped with the
usual topology into the group of linear isometries
of $E$ as before).

A closed subgroup $G$ of ${\rm Iso}(X)$ is a locally compact and
$\sigma$-compact topological group 
and hence it admits a left invariant locally
finite Haar measure $\mu$. In particular, for every $p\in (1,\infty)$
the separable Banach space $L^p(G,\mu)$ of
functions on $G$ which are $p$-integrable with respect to $\mu$
is a Banach module for $G$ 
with respect to the isometric action of $G$ by left translation.
Extending an earlier result for isometry groups of 
proper hyperbolic spaces \cite{H08b} (see also the work of 
Monod-Shalom \cite{MS04}, of 
Mineyev-Monod-Shalom \cite{MMS04}, of Bestvina-Fujiwara \cite{BF08}
and of Caprace-Fujiwara \cite{CF08} 
for closely related results)
we obtain the following non-vanishing result
for second bounded cohomology.

\begin{theorem}\label{thm2} 
Let $G$ be a closed non-elementary subgroup of the isometry
group of a proper ${\rm CAT}(0)$-space $X$
with limit set $\Lambda\subset \partial X$.
If $G$
contains a rank-one element
then $H_{b}^2(G,L^p(G,\mu))\not=\{0\}$ for every $p\in (1,\infty)$.
If $G$ acts transitively on 
the complement of the diagonal in 
$\Lambda\times \Lambda$ 
then we also have $H_{cb}^2(G,L^p(G,\mu))\not=0$.
\end{theorem}

As an application of Theorem \ref{thm2} we
obtain the following super-rigidity theorem.

\begin{coro}
Let $G$ be a connected semi-simple Lie group with
finite center, no compact factors and of rank
at least $2$. Let $\Gamma$ be an irreducible lattice
in $G$, let $X$ be a proper ${\rm CAT}(0)$-space
and let $\rho:\Gamma\to {\rm Iso}(X)$ be a
homomorphism. Let $H<{\rm Iso}(X)$ be the closure of 
$\rho(\Gamma)$. If $H$ is non-elementary and contains
a rank-one element, then $H$ is compact extension of a 
simple
Lie group $L$ of rank one, and up to passing to an open
subgroup of finite index, 
$\rho$ extends to a continuous
homomorphism $G\to L$.
\end{coro}

{\bf Remark:} As in \cite{H08b}, our proof of Theorem \ref{thm2}
also shows the following. Let $G<{\rm Iso}(X)$
be a closed non-elementary subgroup with limit set $\Lambda$ 
which contains a rank-one element.
If $G$ does not act transitively on the complement of
the diagonal in $\Lambda\times \Lambda$ then the second
bounded cohomology group $H_b^2(G,\mathbb{R})$ is infinite dimensional.
However, this was proved by Bestvina and Fujiwara \cite{BF08}.
Moreover, the arguments in \cite{H08b} 
show together with the geometric discussion in Sections 2-5 of this paper
that if $G$ acts transitively on the complement of the
diagonal in $\Lambda\times \Lambda$ then 
$H_{cb}^2(G,\mathbb{R})=0$. Under the additional assumption that
$G$ acts on $X$ cocompactly, this is due to Caprace and Fujiwara \cite{CF08}.

\bigskip

The organization of this note is as follows. In Section 2 we 
collect some geometric properties of proper ${\rm CAT}(0)$-spaces
needed in the sequel. In particular, we discuss contracting geodesics. 
In Section 3 we equip the space of all pairs of endpoints
of uniformly contracting geodesics in $X$ with a family of finite distance
functions parametrized by the points in $X$ which
are equivariant under the natural action 
of the isometry group of $X$.
Section 4 contains the main technical result of the paper.
For a closed non-elementary
subgroup $G$ of ${\rm Iso}(X)$ with limit
set $\Lambda$ which contains a rank one element 
we construct
continuous bounded cocycles for $G$ on a $G$-invariant closed
subspace
of the space of triples of pairwise distinct points in $\Lambda$ 
with values in $L^p(G\times G,
\mu\times \mu)$. This is then used in Section 5 
and Section 6 to show Theorem 2. 
The proof of Theorem 1 and of the corollaries is contained in 
Section 7.

\section{Metric contraction in ${\rm CAT}(0)$-spaces}

In this section we collect some geometric properties
of ${\rm CAT}(0)$-spaces needed in the later
sections. We use
the books \cite{B95,BGS85,BH99} as our main references.

\subsection{Shortest distance projections}

A proper ${\rm CAT}(0)$-space has strong convexity properties
which we summarize in this subsection.

In a complete ${\rm CAT}(0)$-space $X$,
any two points can
be connected by a unique geodesic which varies
continuously with the endpoints. The distance
function is convex: If $\gamma,\zeta:J\to X$ are two
geodesics in $X$ parametrized on the same interval 
$J\subset \mathbb{R}$ then the function 
$t\to d(\gamma(t),\zeta(t))$ is convex.
More generally, we call a function
$f:X\to \mathbb{R}$ \emph{convex} if
for every geodesic $\gamma:J\to \mathbb{R}$ the
function $t\to f(\gamma(t))$ is convex \cite{B95}.

The \emph{visual boundary} $\partial X$ of $X$
is defined to be the space
of all geodesic rays issuing from a fixed point 
$x\in X$ equipped with the
topology of uniform convergence on compact sets.
This definition is independent of the choice of $x$.
We denote the point in $\partial X$ defined by 
a geodesic ray $\gamma:[0,\infty)\to X$ by $\gamma(\infty)$.
We also say that $\gamma$ 
\emph{connects} $x$ to 
$\gamma(\infty)$.

There is another description of the visual boundary of $X$ 
as follows.
Let $C(X)$ be the space of all continuous functions on $X$
endowed with the topology of uniform convergence on
bounded sets. 
Fix a point $y\in X$ and for $x,z\in X$ define
\[b_x(y,z)=d(x,z)-d(x,y).\]
Then we have 
\begin{equation}\label{sym}
b_x(y,z)=-b_x(z,y) \text{ for all } y,z\in X
\end{equation}
and 
\begin{equation}\label{lip}
\vert b_x(y,z)-b_x(y,z^\prime)\vert \leq d(z,z^\prime)
\text{ for all }z,z^\prime\in X\end{equation}
and hence the function 
$b_{x}(y,\cdot):z\to b_x(y,z)$ is 
one-Lipschitz and vanishes at $y$.
Moreover, the function
$b_x(y,\cdot)$ is convex. If $\tilde y\in X$ is 
another basepoint then we have 
\begin{equation}\label{add}
b_x(\tilde y,\cdot)=b_x(y,\cdot)+b_x(\tilde y,y).\end{equation}

The assignment $x\to b_x(y,\cdot)$ is an embedding of $X$ into
$C(X)$. A sequence $\{x_n\}\subset X$ \emph{converges at infinity}
if $d(x_n,y)\to \infty$ and if the functions 
$b_{x_n}(y,\cdot)$ converge in $C(X)$. The visual boundary 
$\partial X$ of $X$ can also be defined as the 
subset of $C(X)$ of all functions 
which are obtained
as limits of functions $b_{x_n}(y,\cdot)$ for
sequences $\{x_n\}\subset X$ 
which converge at infinity. In particular, 
the union $X\cup \partial X$
is naturally a closed subset of $C(X)$
(Chapter II.1 and II.2 of \cite{B95}).
In this way, each $\xi\in \partial X$ corresponds
to a \emph{ Busemann function}  
$b_\xi(y,\cdot)$ at $\xi$ normalized at $y$. If $\gamma:[0,\infty)\to X$
is the geodesic ray which connects $y$ to $\xi$ then this
Busemann function $b_\xi$ satisfies
$b_\xi(y,\gamma(t))=-t$ for all $t\geq 0$.

From now on let $X$ be a \emph{proper} (i.e. complete and 
locally compact) ${\rm CAT}(0)$-space.
Then $X\cup \partial X$ is compact. A subset $C\subset X$ is 
\emph{convex} if for $x,y\in C$ the geodesic connecting
$x$ to $y$ is contained in $C$ as well. 
For every closed convex set $C\subset X$
and every $x\in X$ there is a unique point
$\pi_C(x)\in C$ of smallest distance to $x$ (Proposition II.2.4 of 
\cite{BH99}).
Now let $J\subset \mathbb{R}$ be a closed connected set and
let $\gamma:J\to X$ be a geodesic arc. Then $\gamma(J)\subset X$
is closed and convex and hence there is a 
shortest distance projection $\pi_{\gamma(J)}:X\to \gamma(J)$.
The projection point $\pi_{\gamma(J)}(x)$ of $x\in X$ 
is the unique minimum for the
restriction of the function $b_x(y,\cdot)$ to $\gamma(J)$.
By equality (\ref{add}), 
this does not depend on the choice of the basepoint $y\in X$.
The projection $\pi_{\gamma(J)}:X\to \gamma(J)$ is 
distance non-increasing.

For $\xi\in \partial X$
the function $t\to b_\xi(y,\gamma(t))$
is convex. Let $\overline{\gamma(J)}$ be the closure of $\gamma(J)$
in $X\cup \partial X$. If $b_\xi(y,\cdot)\vert \gamma(J)$ 
assumes a
minimum then we can define $\pi_{\gamma(J)}(\xi)\subset \overline{\gamma(J)}$ 
to be the closure in $\overline{\gamma(J)}$ of the  
connected subset of $\gamma(J)$ of all such minima. 
If $b_\xi(y,\cdot)\vert \gamma(J)$ does not assume a minimum
then by continuity the set $J$ is unbounded and by convexity 
either $\lim_{t\to \infty}b_\xi(y,\gamma(t))=
\inf\{b_\xi(y,\gamma(s))\mid s\in J\}$ 
or $\lim_{t\to -\infty}b_\xi(y,\gamma(t))=
\inf\{b_\xi(y,\gamma(s))\mid s\in J\}$. In the first
case we define 
$\pi_{\gamma(J)}(\xi)=\gamma(\infty)\in \partial X$, and in the
second case we define  
$\pi_{\gamma(J)}(\xi)=\gamma(-\infty)$. Then for every $\xi\in \partial X$
the set  
$\pi_{\gamma(J)}(\xi)$ is a closed connected 
subset of $\overline{\gamma(J)}$ (which may
contain points in both $X$ and $\partial X$).

\subsection{Contracting geodesics}

A ${\rm CAT}(0)$-space may have many totally
geodesic embedded flat subspaces, but it
may also have subsets with hyperbolic behavior.
To give a precise description of such hyperbolic
behavior, Bestvina and Fujiwara introduced
a geometric property for geodesics in 
a ${\rm CAT}(0)$-space (Definition 3.1 of \cite{BF08}) 
which we repeat in the
following definition. For the remainder of this note,
geodesics are always defined on closed connected
subsets of $\mathbb{R}$.

\begin{definition}\label{contracting}
A geodesic arc $\gamma:J\to X$ is 
\emph{$B$-contracting}
for some $B>0$ if for every closed 
metric ball $K$ in $X$ which is disjoint
from $\gamma(J)$ the diameter of the 
projection $\pi_{\gamma(J)}(K)$ does not exceed $B$. 
\end{definition}

We call a geodesic \emph{contracting} if it is $B$-contracting
for some $B>0$. Lemma 3.3 of \cite{H08d} relates
$B$-contraction for a geodesic $\gamma$ to the diameter of the projections
$\pi_{\gamma(\mathbb{R})}(\xi)$ where $\xi\in \partial X$.

\begin{lemma}\label{diamproj}
Let $\gamma:\mathbb{R}\to X$ be a $B$-contracting geodesic.
Then for every $\xi\in \partial X-\{\gamma(-\infty),\gamma(\infty)\}$
the projection $\pi_{\gamma(\mathbb{R})}(\xi)$ is a 
compact subset of $\gamma(\mathbb{R})$ of diameter at most $6B+4$.
\end{lemma}

Lemma 3.2 and 3.5 of \cite{BF08} show that 
a connected subarc of a contracting geodesic is contracting and that
a triangle containing a $B$-contracting
geodesic as one of its sides is
uniformly thin.

\begin{lemma}\label{thintriangle} 
\begin{enumerate}
\item Let $\gamma:J\to X$ be a $B$-contracting
geodesic. Then for every closed connected subset $I\subset J$,
the subarc $\gamma(I)$ of $\gamma$ is $B+3$-contracting.
\item
Let $\gamma:J\to X$ be a $B$-contracting geodesic.
Then for $x\in X$ and for every
$t\in J$ the geodesic connecting $x$ to $\gamma(t)$
passes through the $3B+1$-neighborhood of $\pi_{\gamma(J)}(x)$. 
\end{enumerate}
\end{lemma}

Note that by convexity of the distance function,
if $\zeta_i:[a_i,b_i]\to X$ $(i=1,2)$ are two geodesic
segments such that $d(\zeta_1(a_1),\zeta_2(a_2))\leq R,
d(\zeta_1(b_1),\zeta_2(b_2))\leq R$ then the
\emph{Hausdorff distance} between the subsets
$\zeta_1[a_1,b_1],\zeta_2[a_2,b_2]$ of $X$ is at most $R$.
Here the Hausdorff distance between closed 
(not necessarily compact)
subsets $A,B$ of $X$ is the infimum of all numbers $R>0$
such that $A$ is contained in the $R$-neighborhood of $B$
and $B$ is contained in the $R$-neighborhood of $A$.
(This number may be infinite).

A \emph{visibility point} is a point $\xi\in \partial X$
with the property that any $\eta\in \partial X-\{\xi\}$ can
be connected to $\xi$ by a geodesic line.
By Lemma 3.5 of \cite{H08d},  
the endpoint of a contracting geodesic ray is a visibility point.
Geodesic rays which abut at an endpoint of a 
contracting geodesic ray are themselves contracting.

\begin{lemma}\label{localray}
For every $B>0$ there is a number $C=C(B)>B$ with the
following property.
Let $\gamma:[0,\infty)\to X$ be a $B$-contracting ray and
let $\xi\in \partial X-\gamma(\infty)$. Then every
geodesic $\zeta$ connecting 
$\xi=\zeta(-\infty)$ to $\gamma(\infty)=\zeta(\infty)$ 
passes through the $9B+6$-neighborhood of 
every point $x\in \pi_{\gamma[0,\infty)}(\xi)$.
If $t\in\mathbb{R}$ is such that
$d(\zeta(t),x)\leq 9B+6$ then the
geodesic ray $\zeta[t,\infty)$ is $C$-contracting.
\end{lemma}
\begin{proof}
Let $\gamma:[0,\infty)\to X$ be a $B$-contracting geodesic ray
and let $\xi\in \partial X-\gamma(\infty)$. 
Assume that $\gamma(s)\in \pi_{\gamma[0,\infty)}(\xi)$.
By Lemma \ref{diamproj},  
the projection $\pi_{\gamma[0,\infty)}(\xi)$ 
is contained in 
$\gamma[s-6B-4,s+6B+4]$.
Let $\zeta:\mathbb{R}\to X$ be a geodesic connecting $\xi$ 
to $\gamma(\infty)$. 
Lemma 2.2 of \cite{H08d} shows
that for sufficiently large $t$ we have
\[\pi_{\gamma[0,\infty)}(\zeta(-t))\in 
\gamma[s-6B-5,s+6B+5].\]
Thus by Lemma \ref{thintriangle}, the geodesic 
ray $\zeta[-t,\infty)$ connecting $\zeta(-t)$ 
to $\gamma(\infty)$ 
passes through the $9B+6$-neighborhood of $\gamma(s)$. 
Since $\zeta$ was an arbitrary geodesic connecting
$\xi$ to $\gamma(\infty)$, the first part of the lemma follows.

From this and Lemma 3.8 of \cite{BF08}, the
second part of the lemma is immediate as well.
Namely, let again $\zeta $ be a geodesic connecting
$\xi$ to $\gamma(\infty)$ and assume that
$\zeta$ is parametrized in such a way
that $d(\zeta(0),\pi_{\gamma[0,\infty)}(\xi))\leq 9B+6$. 
The geodesic ray $\zeta[0,\infty)$ is a locally uniform limit as
$t\to \infty$ of the geodesics
$\zeta_t$ connecting $\zeta(0)$ to 
$\gamma(t)$. By Lemma 3.8 of \cite{BF08},
there is a number $C>0$ only depending on $B$ such that 
each of the geodesics $\zeta_t$ is $C$-contracting.
Now Lemma 3.3 of \cite{H08d} shows that
a limit of a sequence of $C$-contracting
geodesics is $C$-contracting from which the 
lemma follows.
\end{proof}

The next observation is an extension of Lemma \ref{thintriangle}.
For its formulation, define an \emph{ideal geodesic triangle}
to consist of three biinfinite geodesics $\gamma_1,\gamma_2,\gamma_3$
with $\gamma_i(\infty)=\gamma_{i+1}(-\infty)$ (where indices are taken
modulo three). The points $\gamma_i(\infty)$ $(i=1,2,3)$ are called
the vertices of the ideal geodesic triangle.
If $a,b\in \partial X$ are visibility
points then for every $\xi\in \partial X-\{a,b\}$ there is an ideal
geodesic triangle with vertices $a,b,\xi$. Note that such a triangle
need not be unique.

\begin{lemma}\label{thintriangle2}
Let $B>0$ and
let $\gamma:\mathbb{R}\to X$ be a $B$-contracting geodesic.
Then for every ideal geodesic
triangle $T$ with side $\gamma$ there
is a point $x\in X$ whose distance to each of the sides of $T$ does not
exceed $9B+6$. The diameter of the set of all such
points does not exceed $54B+36$.
\end{lemma}
\begin{proof}
Let $\gamma:\mathbb{R}\to X$ be a $B$-contracting geodesic and let 
$T$ be an ideal geodesic triangle with side $\gamma$ and 
vertex $\xi\in \partial X-\{\gamma(\infty),\gamma(-\infty)\}$
opposite to $\gamma$. Assume that $\gamma$ is parametrized
in such a way that  
$\gamma(0)\in \pi_{\gamma(\mathbb{R})}(\xi)$. 
Let $c:[0,\infty)\to X$ be the geodesic ray connecting
$c(0)=\gamma(0)$ to $\xi$. By Lemma \ref{localray}
and by ${\rm CAT}(0)$-comparison, the side $\alpha$ of $T$ 
connecting $\xi$ to $\gamma(\infty)$ is contained
in the $9B+6$-tubular neighborhood of $c[0,\infty)\cup
\gamma[0,\infty)$, and the side $\beta$ of $T$ connecting
$\xi$ to $\gamma(-\infty)$ is contained in the
$9B+6$-tubular neighborhood of $c[0,\infty)\cup
\gamma(-\infty,0]$. The distance
between $\gamma(0)=c(0)$ and every side of $T$ does not
exceed $9B+6$.

Since the projection $\pi_{\gamma(\mathbb{R})}$ is distance
non-increasing and since $\pi_{\gamma(\mathbb{R})}(c[0,\infty))=
\gamma(0)$ (see \cite{BH99} and
the proof of Lemma 3.5 in \cite{H08d}),
if as before $\alpha,\beta$ are the sides of $T$ connecting
$\xi$ to $\gamma(\infty),\gamma(-\infty)$, respectively, then
\[\pi_{\gamma(\mathbb{R})}(\alpha)\subset \gamma[-9B-6,\infty)
\text{ and }
\pi_{\gamma(\mathbb{R})}(\beta)\subset \gamma(-\infty,9B+6].\]
Now if $x\in X$ is such that 
the distance between $x$ and each side of $T$ is at most
$9B+6$ then using again that
$\pi_{\gamma(\mathbb{R})}$ is distance
non-increasing we conclude that $\pi_{\gamma(\mathbb{R})}(x)\in
\gamma[-18B-12,18B+12]$. But $d(x,\gamma(\mathbb{R}))\leq 9B+6$ and hence
$d(x,\gamma(0))\leq 27B+18$.
This completes the proof of the lemma.
\end{proof}

\subsection{Isometries}

For an isometry $g$ of $X$ define the displacement function
$d_g$ of $g$ to be the function $x\to d_g(x)=d(x,gx)$.
An isometry $g$ of $X$ 
is called \emph{semisimple} if $d_g$ assumes
a minimum in $X$. If $g$ is semisimple and ${\rm min}\,d_g=0$
then $g$ is called \emph{elliptic}. 
Thus an isometry is elliptic if and only
if it fixes at least one point in $X$. 
A semisimple isometry
$g$ with ${\rm min}\,d_g>0$ is called \emph{axial}.
By Proposition 3.3 of \cite{B95}, an isometry $g$ of $X$ is
axial if and only if there is a geodesic
$\gamma:\mathbb{R}\to X$ such that
$g\gamma(t)=\gamma(t+\tau)$ for every $t\in \mathbb{R}$ where
$\tau=\min\,d_g>0$ is the \emph{translation length} of $g$.
Such a geodesic is called an \emph{oriented axis} for $g$.
Note that the geodesic $t\to \gamma(-t)$ is an oriented 
axis for $g^{-1}$. The endpoint $\gamma(\infty)$ of $\gamma$
is a fixed point for the action of $g$ on $\partial X$
which is called the \emph{attracting fixed point}. 
The closed convex set $A\subset X$ of all points for which the displacement
function of $g$ is minimal is 
isometric to a product space $C\times \mathbb{R}$.
For each
$x\in C$ the set $\{x\}\times \mathbb{R}$ is an axis of $g$.

Bestvina and Fujiwara
introduced the following notion to identify 
isometries of a ${\rm CAT}(0)$-space with geometric
properties similar to the properties of isometries
in a hyperbolic geodesic metric space
(Definition 5.1 of \cite{BF08}).

\begin{definition}\label{rankone}
An isometry $g\in {\rm Iso}(X)$ is called 
\emph{$B$-rank-one} for some $B>0$ 
if $g$ is axial and admits a $B$-contracting axis.
\end{definition}
We call an isometry $g$ \emph{rank-one} if $g$ is $B$-rank-one for some
$B>0$.

The following statement is Theorem 5.4 of \cite{BF08}.

\begin{proposition}\label{flatrank}
An axial isometry of $X$ with axis $\gamma$ is rank-one
if and only if $\gamma$ does not bound a flat half-plane.
\end{proposition}

Let $G<{\rm Iso}(X)$ be a subgroup of the isometry group of $X$.
The \emph{limit set} $\Lambda$ of $G$ is the set
of accumulation points in $\partial X$ of one (and hence every)
orbit of the action of $G$ on $X$. The
limit set is a compact non-empty $G$-invariant subset of $\partial X$. 
Call $G$ \emph{non-elementary} if its limit set contains
at least three points and if moreover $G$ does not fix
globally a point in $\partial X$.

A compact space is \emph{perfect} if it does not
have isolated points. The action of a group $G$ on 
a topological space $\Lambda$ is called \emph{minimal} if 
every orbit is dense.
A homeomorphism $g$ of a space $\Lambda$ is said
to act with \emph{north-south dynamics} if there are
two fixed points $a\not=b\in \Lambda$ for the action of $g$ such that
for every neighborhood $U$ of $a$, $V$ of $b$ there is some 
$k>0$ such that $g^k(\Lambda-V)\subset U$ and 
$g^{-k}(\Lambda-U)\subset V$.
The point $a$ is called the \emph{attracting fixed point} for $g$,
and $b$ is the repelling fixed point. The following is shown in 
\cite{H08d} (see also \cite{B95} for a similar discussion).

\begin{lemma}\label{northsouth}
Let $G<{\rm Iso}(X)$ be a non-elementary group which
contains a rank-one element. Then the limit set $\Lambda$ of $G$ is
perfect, and it is the smallest closed $G$-invariant subset of
$\partial X$. The action of $G$ on $\Lambda$ is minimal. 
An element $g\in G$ is rank-one
if and only if $g$ acts on $\partial X$ with
north-south dynamics.
\end{lemma}

For every proper metric space $X$, the isometry group ${\rm
Iso}(X)$ of $X$ can be equipped with a natural locally compact
$\sigma$-compact metrizable topology, the so-called \emph{compact
open topology}. With respect to this topology, a sequence
$(g_i)\subset {\rm Iso}(X)$ converges to some isometry $g$ if and
only if $g_i\to g$ uniformly on compact subsets of $X$.
In this topology, a
closed subset $A\subset {\rm Iso}(X)$ is compact if and only if
there is a compact subset $K$ of $X$ such that $gK\cap
K\not=\emptyset$ for every $g\in A$. In particular, the action of
${\rm Iso}(X)$ on $X$ is proper.
In the
sequel we always equip subgroups of ${\rm Iso}(X)$
with the compact open topology.

Denote by $\Delta$ the diagonal in 
$\partial X\times \partial X$. 
We have (Lemma 6.1 of \cite{H08b}).

\begin{lemma}\label{closed}
Let $G<{\rm Iso}(X)$ be a closed subgroup with limit set $\Lambda$.
Let
$(a,b)\in \Lambda\times \Lambda-\Delta$ be the pair of fixed points of a
rank one element of $G$. 
Then the $G$-orbit of $(a,b)$ is a closed
subset of $\Lambda\times \Lambda-\Delta$.
\end{lemma}

The following technical observation is useful
in Section 6.

\begin{lemma}\label{compact}
Let $G<{\rm Iso}(X)$ be a closed non-elementary group 
with limit set $\Lambda$. If $G$ contains
a rank-one element $g\in G$ with 
fixed points $a\not= b\in \Lambda$ and   
if $G$ does not act transitively on the complement of the diagonal
in $\Lambda\times\Lambda$ then
there is some $h\in G$ such that the stabilizer in $G$
of the pair of points $(b,hb)\in \Lambda\times\Lambda$
is compact.  
\end{lemma}
\begin{proof}
Let $g\in G$ be 
a rank-one element with
attracting fixed point $a\in \Lambda$,
repelling fixed point $b\in \Lambda$. Let
$\gamma:\mathbb{R}\to X$ be an axis for $g$ connecting
$b$ to $a$. Then $\gamma$ is $B$-contracting for some $B>0$.

Let $h\in G$ be such that $hb\not=b$. 
The rays $\gamma(-\infty,0]$ and $h(\gamma(-\infty,0])$ are
$B$-contracting, with endpoints 
$b,hb\in \Lambda$. Now a biinfinite geodesic $\xi:\mathbb{R}\to X$
with the property that there are numbers 
$-\infty<s<t<\infty, C>0$ such that the rays 
$\xi(-\infty,s],\xi[t,\infty)$ are $C$-contracting is
$C^\prime$-contracting for a number $C^\prime >C$
only depending on $C$ and on $[s,t]$. Therefore
Lemma \ref{localray} implies that
there is a number $B_0>B$ such that each geodesic
connecting $b$ to $hb$ is $B_0$-contracting. 
As a consequence, the set $A\subset X$ of all
points which are contained in 
a geodesic connecting $b$ to $hb$ is 
closed and convex and isometric to 
$K_0\times \mathbb{R}$ for a 
compact convex subset $K_0$ of $X$.

An isometry of $X$ which stabilizes the pair of 
points $(b,hb)$ preserves the closed convex set
$K_0\times \mathbb{R}$. In particular, each such isometry $u$ 
is semi-simple. Moreover, if $u$ is not elliptic then
$u$ is rank-one. Since $G$ is a closed subgroup of ${\rm Iso}(X)$,
this implies that either the stabilizer of $(b,hb)$ in $G$ is compact
or it contains a rank one element.

Now assume that there is no
$h\in G$ such that the stabilizer of 
$(b,hb)$ in $G$ is compact. Then each such stabilizer contains
a rank-one element. The stabilizer $G_b$ of $b$
in $G$ is a closed subgroup of $G$.

Let $h\in G$ with $hb\not\in\{a,b\}$ and 
let again $\gamma$ be an oriented axis for $g$
connecting $b$ to $a$. Assume that $\gamma(0)\in
\pi_{\gamma(\mathbb{R})}(hb)$. Since $\gamma$ is $B$-contracting,
by Lemma \ref{localray} and convexity the ray
$\gamma(-\infty,0]$ is contained in 
the $9B+6$-neighborhood
of every geodesic connecting $b$ to $hb$.
By assumption and the above discussion, there 
is a rank-one element $u\in G_b$ with attracting
fixed point $hb$ and repelling fixed point $b$.
Since $u$ acts with north-south dynamics on $\Lambda$, we have
$u^ia\to hb$ $(i\to \infty)$ and hence also
$u^ig^{-k}a\to hb$ for all $k$.
Let $\tau>0$ be the translation length of $g$ and
let $K$ be the
closed $18B+12+2\tau$-neighborhood of $\gamma(0)$. 
By the choice of the set $K$, for every $i>0$ there is
some $k=k(i)>0$ such that $u^ig^{-k}\gamma(0)\in K$. 
Since $G$ is a closed subgroup of ${\rm Iso}(X)$, up to
passing to a subsequence the sequence
$\{u^ig^{-k}\}\subset G_b$ converges to an element $v\in G_b$
with $va=hb$. 

By the discussion in the previous paragraph, 
for every $x\in Gb-\{b\}$ there is some $v\in G_b$ with
$va=x$. This implies that 
the image of $(a,b)$ under the action of the
group $G$ is dense in $\Lambda\times \Lambda-\Delta$. 
Namely, by Lemma \ref{northsouth}, the $G$-orbit of $b$ is 
dense in $\Lambda$. Thus it suffices to show that for every
$u\in G$ and every $x\in Gb-\{b,ub\}$ there is
some $v\in G$ with $v(a,b)=(x,ub)$. For this let 
$y=u^{-1}x$. Then $y\not=b$ and hence there is some
$w\in G_b$ with $w(a)=y$. Then the isometry
$uw$ satisfies $uw(b)=u(b)$, $uw(a)=x$.

Now by Lemma \ref{closed}, the $G$-orbit of $(a,b)$ is 
closed in $\Lambda\times \Lambda-\Delta$ and hence together
with the above we conclude that $G$ acts transitively
on the complement of the diagonal in $\Lambda\times \Lambda$. 
The lemma follows.
\end{proof}

A free group with two generators is hyperbolic
in the sense of Gromov \cite{GH}. In particular,
it admits a Gromov boundary which can be viewed
as a compactification of the group. The following result 
is contained in \cite{BF08} (see also Proposition 5.8
of \cite{H08d} and \cite{CF08,BF02}).

\begin{lemma}\label{free}
Let $G<{\rm Iso}(X)$ be a closed non-elementary group which 
contains a rank-one element. Let $\Lambda\subset \partial X$
be the limit set of $G$. 
If $G$ does not act transitively on $\Lambda\times \Lambda-\Delta$
then $G$ contains a free subgroup $\Gamma$ with two generators
and the following properties.
\begin{enumerate}
\item Every element $e\not=g\in \Gamma$ is rank-one.
\item There is a $\Gamma$-equivariant embedding of
the Gromov boundary of $\Gamma$ into $\Lambda$.
\item  There are infinitely many
elements $u_i\in \Gamma$ $(i>0)$ with 
fixed points $a_i,b_i$
such that for all $i$ the $G$-orbit of $(a_i,b_i)
\in \Lambda\times \Lambda-\Delta$ is distinct from the orbit of
$(b_j,a_j) (j>0)$ or $(a_j,b_j)(j\not=i)$.
\end{enumerate}
\end{lemma}

\section{The space of $B$-contracting geodesics}

In the previous section we introduced for some $B>0$ 
a $B$-contracting geodesic in a proper ${\rm CAT}(0)$-space $X$.
In this section we consider in more
detail the space of all such geodesics in $X$.

The main idea is as follows: Even though the 
geometry of a ${\rm CAT}(0)$-space
$X$ may be very different from the geometry of 
a hyperbolic geodesic metric space, if $X$ admits
$B$-contracting geodesics then by Lemma \ref{thintriangle2},
these geodesics have the same global geometric
properties as geodesics in a $\delta$-hyperbolic
geodesic metric space where $\delta >0$ only
depends on $B$. As a consequence, given a fixed
point $x\in X$, we can describe the position of 
two such
geodesics $\gamma,\zeta$ relative to each other as seen
from $x$ by introducing a metric quantity which can be
thought of being equivalent to the (oriented) sum of the
Gromov distances at $x$ of their endpoints in the
case that the space $X$ is hyperbolic.

We continue
to use the assumptions and notations from Section 2. 
In the remainder of this section, a geodesic in $X$
is always defined on a closed connected subset $J$ of $\mathbb{R}$. 
For some $B>0$ denote by 
${\cal A}(B)\subset \partial X\times \partial X-\Delta$  
the set of all pairs of points in $\partial X$
which are connected
by a $B$-contracting geodesic. We have.

\begin{lemma}\label{abclosed}
${\cal A}(B)$ is a closed subset of $\partial X\times \partial X-\Delta$.
\end{lemma}
\begin{proof} Let $\{(\xi_i,\eta_i)\}\subset {\cal A}(B)$ 
be a sequence which converges in $\partial X\times \partial X-\Delta$
to a point $(\xi,\eta)$. For each $i$  
let $\gamma_i$ be a $B$-contracting geodesic connecting $\xi_i$ to $\eta_i$.
We first claim that the geodesics $\gamma_i$ pass through a 
fixed compact subset of $X$.

Namely, choose a point $x\in X$ and let 
$x_i=\pi_{\gamma_i(\mathbb{R})}(x)$. 
If the geodesics $\gamma_i$ do not pass through a fixed
compact subset of $X$ then we have 
$d(x_i,x)\to \infty$. 
Since $X\cup \partial X$ is compact, 
after passing to a subsequence we may assume that $x_i\to \alpha\in \partial X$
as $i\to \infty$. 
On the other hand, the geodesic $\gamma_i$ is
$B$-contracting and therefore by Lemma \ref{thintriangle}
the geodesics connecting $x$ to 
$\xi_i=\gamma_i(-\infty),\eta_i=\gamma_i(\infty)$
both pass through the $3B+1$-neighborhood of $x_i$.
By ${\rm CAT}(0)$-comparison, this implies that 
$\xi_i\to \alpha,\eta_i\to \alpha$ which
contradicts the assumption that $\xi_i\to \xi,\eta_i\to 
\eta\not=\xi$.  

Thus the geodesics $\gamma_i$ pass through a fixed compact
subset of $X$ and therefore after passing to a subsequence
we may assume that $\gamma_i\to\gamma$ locally uniformly
where $\gamma$ is a geodesic connecting $\xi$ to $\eta$.
The limit geodesic is $B$-contracting
by Lemma 3.6 of \cite{H08d}.
\end{proof}

For a number $B>0$, 
a point $x\in X$ and an ordered pair
$(\zeta_1:J_1\to X,\zeta_2:J_2\to X)$ of oriented geodesics in $X$ 
which share at most one endpoint in $\partial X$
define a number $\tau_B(x,\zeta_1,\zeta_2)\geq 0$ 
as follows.

By convexity of the distance function, there are (perhaps empty)
closed connected subsets $[a_1,b_1]\subset J_1,
[a_2,b_2]\subset J_2$ such that
\[[a_i,b_i]=\{t\mid d(\zeta_i(t),\zeta_{i+1}(J_{i+1}))\leq 6B+2\}.\]
(Here $i=1,2$ and indices are taken modulo two.
If $\zeta_1,\zeta_2$ have a common endpoint in $\partial X$ then 
one of the numbers $a_1,b_1$ and one of the numbers
$a_2,b_2$ may be infinite.) 

If $[a_i,b_i]\not=\emptyset$ then 
let $s_i,t_i\in J_i\cup\{\pm \infty\}$ be such that
\[\pi_{\zeta_i(J_i)}(\zeta_{i+1}(a_{i+1}))=\zeta_i(s_i)
\text{ and }\pi_{\zeta_i(J_i)}(\zeta_{i+1}(b_{i+1}))=\zeta_i(t_i)\] 
($i=1,2$ and indices are taken modulo two)
and let $x_i=\pi_{\zeta_i(J_i)}(x)$
$(i=1,2)$. 

If $s_i<t_i$ and if $x_i\in \zeta_i[a_i,b_i]$ for $i=1,2$ then 
define 
\[\tau_B(x,\zeta_1,\zeta_2)=
\min\{d(x_i,\zeta_i(a_i)),d(x_i,\zeta_i(b_i))\mid i=1,2\}.\]
In all other cases define $\tau_B(x,\zeta_1,\zeta_2)=0$.
Note that $\tau_B(x,\zeta_1,\zeta_2)$ 
depends on the orientation of $\zeta_1,\zeta_2$ but not
on the parametrization of $\zeta_1,\zeta_2$ defining a fixed orientation.

We collect some first easy properties of the 
function $\tau_B$.

\begin{lemma}\label{taubfirst}
For any two geodesics 
$\zeta_1,\zeta_2$ in $X$ and any
$x\in X$ the following holds true.
\begin{enumerate}
\item 
$\tau_B(x,\zeta_1,\zeta_2)=\tau_B(x,\zeta_2,\zeta_1)$.
\item If $\hat \zeta_i$ equals the geodesic obtained from 
$\zeta_i$ by reversing the orientation then 
$\tau_B(x,\hat\zeta_1,\hat \zeta_2)=\tau_B(x,\zeta_1,\zeta_2)$.
\item $\tau_B(x,\zeta_1,\zeta_2)\leq \tau_B(y,\zeta_1,\zeta_2)+d(x,y)$
for all $x,y\in X$.
\end{enumerate}
\end{lemma}
\begin{proof}
The first and the second property in the lemma is obvious from the
definition. To show the third property simply
note that for a geodesic $\zeta:J\to X$ 
the projection $\pi_{\zeta(J)}$ is
distance non-increasing. 
\end{proof}

Moreover we observe.

\begin{lemma}\label{tauproj}
Let $\zeta_i:J_i\to X$ 
be $B$-contracting geodesics $(i=1,2)$ such that
$\tau_B(x,\zeta_1,\zeta_2)>0$. Then we have
\[d(\pi_{\zeta_1(J_1)}(x),
\pi_{\zeta_2(J_2)}(x))\leq 24B+8.\]
\end{lemma}
\begin{proof} If $\tau_B(x,\zeta_1,\zeta_2)>0$ and if
(after reparametrization) we have 
$\pi_{\zeta_i(J_i)}(x)=\zeta_i(0)$ $(i=1,2)$ then 
by definition of the function $\tau_B$ there is a point
on the geodesic $\zeta_2$ whose distance to 
$\zeta_1(0)$ does not exceed $6B+2$. 
This shows that
\[d(x,\zeta_2(J_2))\leq d(x,\zeta_1(J_1)) +6B+2.\]
By symmetry we conclude that
\begin{equation}\label{abs}
\vert d(x,\zeta_1(J_1))-d(x,\zeta_2(J_2))\vert \leq 6B+2.
\end{equation}
Thus if $t\in J_2$ is such that $d(\zeta_1(0),\zeta_2(t))\leq 6B+2$
then 
\begin{equation}\label{below}
d(x,\zeta_2(t))\leq d(x,\zeta_2(J_2))+12B+4.
\end{equation} 

On the
other hand, by Lemma \ref{thintriangle}, 
the geodesic connecting $x$ to $\zeta_2(t)$ passes
through the $3B+1$-neighborhood of $\zeta_2(0)$ and hence
\begin{equation}\label{above}
d(x,\zeta_2(t))\geq d(x,\zeta_2(J_2))+\vert t\vert
-6B-2.\end{equation}
The two inequalities (\ref{below}) and (\ref{above})
together show that $\vert t\vert \leq 18B+6$
and therefore $d(\pi_{\zeta_1(J_1)}(x),
\pi_{\zeta_2(J_2)}(x))\leq 24B+8$ as claimed.
\end{proof}

We also have.

\begin{lemma}\label{decay} 
Let $\zeta_i:[0,\infty)\to X$ $(i=1,2)$ 
be two geodesic rays
with the same endpoint $\zeta_1(\infty)=\zeta_2(\infty)$.
Let $s\in [1,\infty)$ be such that
\[p=\tau_{B}(\zeta_1(s),\zeta_1,\zeta_2)\geq 1.\] Then 
$\tau_{B}(\zeta_1(s+t),\zeta_1,\zeta_2)\geq p+t-12B-4$
for all $t\geq 0$.
\end{lemma}
\begin{proof}
Let $\zeta_i:[0,\infty)\to X$ be geodesic rays in $X$ $(i=1,2)$ 
with $\zeta_1(\infty)=\zeta_2(\infty)$.
Let $s\in [0,\infty)$ be such that
$\tau_B(\zeta_1(s),\zeta_1,\zeta_2)\geq 1$.
Then the geodesic ray $\zeta_2[0,\infty)$ passes through the
$6B+2$-neighborhood of $\zeta_1(s)$. 
If $s^\prime\in[0,\infty)$
is such that $d(\zeta_1(s),\zeta_2(s^\prime))\leq 6B+2$ then
by convexity of the distance function we have
\[d(\zeta_1(s+t),\zeta_2(s^\prime+t))\leq 6B+2\text{ for all }t\geq 0.\]
Now let $t\geq 0$ and let 
$\sigma\in \mathbb{R}$ be such that
$\pi_{\zeta_2[0,\infty)}(\zeta_1(s+t))=\zeta_2(\sigma)$.
Then $d(\zeta_1(s+t),\zeta_2(\sigma))\leq 6B+2$ and 
hence the triangle inequality shows that
$\sigma\in [s^\prime+t-12B-4,s^\prime+t+12B+4]$.
From this and the definition of the function 
$\tau_{B}$ the lemma follows.
\end{proof}

The next observation is the analog of the 
familiar ultrametric inequality for Gromov
products in hyperbolic spaces.

\begin{lemma}\label{tauestimate}
There is a number $L>0$ such that for 
every $B>0$ and for all $B$-contracting
geodesics $\zeta_i:J_i\to X$ $(i=1,2)$ we have 
\[\tau_B(x,\zeta_1,\zeta_3)\geq 
\min\{\tau_B(x,\zeta_1,\zeta_2),\tau_B(x,\zeta_2,\zeta_3)\}-LB.\]
\end{lemma}
\begin{proof}
Let $\zeta_i:J_i\to X$ be $B$-contracting geodesics and let
$x\in X$.
Number the geodesics $\zeta_i$ 
in such a way that $\tau_B(x,\zeta_1,\zeta_3)=
\min\{\tau_B(x,\zeta_i,\zeta_{i+1})\mid
i=1,2,3\}$. 
Assume also without loss of
generality that $r_1=\tau_B(x,\zeta_1,\zeta_2)\leq 
r_2=\tau_B(x,\zeta_2,\zeta_3)$. 
If $r_1=0$ then there is nothing to show. So assume that
$r_1 >0$. 
By Lemma \ref{tauproj}
we then have
\begin{equation}\label{projdist}
d(\pi_{\zeta_{2}(J_2)}(x),\pi_{\zeta_{j}(J_j)}(x))
\leq 24B+8 \,(j=1,3)
\end{equation}
and hence 
$d(\pi_{\zeta_1(J_1)}(x),\pi_{\zeta_3(J_3)}(x))\leq 48B+16.$
Since the lemma is only significant if $r_1$ is large, 
we successively increase a lower bound for $r_1$
by a controlled amount in the course
of the proof 
so that all the geometric estimates are meaningful
without explicit mentioning.

For simplicity parametrize 
the geodesics $\zeta_i$ in such a way
that $\pi_{\zeta_i(J_i)}(x)=\zeta_i(0)$ $(i=1,2,3)$. By definition
of the function $\tau_B$ there is a number $t_2\geq 0$ such that
$d(\zeta_1(r_1),\zeta_2(t_2))\leq 6B +2$. By the distance
estimate (\ref{projdist}), we have 
\[t_2 =d(\zeta_2(t_2),\pi_{\zeta_2(J_2)}(x))
\in [r_1-30B-10,r_1+30B+10]\]
and hence $d(\zeta_1(r_1),\zeta_2(r_1))\leq 36B+12$. 
Now $r_2\geq r_1$ by
assumption and therefore using once more the definition of the 
function $\tau_B$ we have $d(\zeta_2(r_1),\zeta_3(J_3))\leq 6B+2$.
Thus if we write $R_0=42B +18$ then we have
$d(\zeta_1(r_1),\zeta_3(J_3))\leq R_0$ and
similarly $d(\zeta_1(-r_1),\zeta_3(J_3))\leq R_0$. 
Since $d(\pi_{\zeta_1(J_1)}(x),\pi_{\zeta_3(J_3)}(x))\leq 48B+16$,
by the estimate (\ref{projdist}) above 
we conclude that for 
$R_1=48B+16+R_0$ there 
are numbers 
$s_3,t_3\geq r_1-R_1$ such that $[-s_3,t_3]\subset J_3$
and that
\begin{equation}\label{r0}d(\zeta_1(-r_1),\zeta_3(-s_3))\leq R_0
\text{ and }d(\zeta_1(r_1),\zeta_3(t_3))\leq R_0.\end{equation}

By assumption,
$\zeta_1$ and $\zeta_3$ are $B$-contracting.
Let $\rho:[0,b]\to X$ be the geodesic connecting 
$\zeta_1(-r_1)=\rho(0)$ to $\zeta_3(t_3)=\rho(b)$. 
Let $z=\pi_{\zeta_1(J_1)}(\zeta_3(t_3))$. 
By the estimate (\ref{r0}) and the triangle 
inequality, the distance
between $z$ and $\zeta_3(t_3)$ is at most $R_0$,
and the distance
between $z$ and $\zeta_1(r_1)$ is bounded from above by
$2R_0$.

Since $\zeta_1$ is $B$-contracting, by
Lemma \ref{thintriangle} and the remark thereafter,
there is a number $T\leq b$
such that the Hausdorff distance between 
the subarc of 
$\zeta_1$ connecting $\zeta_1(-r_1)$ to $z$ 
and the arc $\rho[0,T]$ is at most $3B+1$.
Moreover, we can choose $T$ in such a way that
$T\geq b-R_0$.

Similarly, 
since $\zeta_3$ is $B$-contracting, if
$w=\pi_{\zeta_3(J_3)}(\zeta_1(-r_1))$ then 
the distance between $w$ and $\zeta_1(-r_1)$ is
at most $R_0$. There is a number
$S\leq R_0$ such that the Hausdorff distance
between $\rho[S,T]$ and the subarc
of $\zeta_3$ connecting $w$ to $\zeta_3(t_3)$ is
at most $3B+1$.

As a consequence, there are two subarcs $\zeta_1^\prime$
of $\zeta_1$, $\zeta_2^\prime$ of $\zeta_2$ whose
Hausdorff distance to the geodesic arc 
$\rho[S,T]$ is at most
$3B+1$. Hence the Hausdorff distance between
$\zeta_1^\prime$ and $\zeta_2^\prime$ is at most $6B+2$.

If $r_1$ is sufficiently large depending
on $B$ then 
there is a number $L>0$ and there is a 
subarc $\zeta_1^\prime,\zeta_3^\prime$
of $\zeta_1(J_1),\zeta_3(J_3)$ with the following property.
The arc $\zeta_i^\prime$ 
contains $x_1,x_3$ as interior points, and the distance
of $x_1,x_3$ to the endpoints of $\zeta_1^\prime,
\zeta_2^\prime$ is 
at least $r_1-LB$. Moreover, the Hausdorff
distance in $X$ between $\zeta_1^\prime,\zeta_3^\prime$ 
is smaller than $6B+2$. This shows that
\[\tau_B(x,\zeta_1,\zeta_3)\geq
r_1-LB\geq \tau_B(x,\zeta_1,\zeta_2)-LB\]  
which completes the proof of the lemma.
\end{proof}

For distinct pairs of points $(\xi_1,\eta_1),
(\xi_2,\eta_2)\in {\cal A}(B)$ define
\[\tau_B(x,(\xi_1,\eta_1),(\xi_2,\eta_2))
\geq 0\] to be the infimum of the numbers
$\tau_B(x,\zeta_1,\zeta_2)$ over all $B$-contracting geodesics 
$\zeta_i$ connecting $\xi_i$ to $\eta_i$ $(i=1,2)$.
Clearly we have \[\tau_B(x,\alpha_1,\alpha_2)=
\tau_B(x,\alpha_2,\alpha_1)\text{ for all }x\in X,
\alpha_1,\alpha_2\in {\cal A}(B).\]
Moreover, by Lemma \ref{tauestimate}, there is a number
$L>0$ such that for 
all $\alpha_1,\alpha_2,\alpha_3\in {\cal A}(B)$ 
we have
\[\tau_B(x,\alpha_1,\alpha_3)\geq
\min\{\tau_B(x,\alpha_1,\alpha_2),\tau_B(x,\alpha_2,\alpha_3)\}-LB.\]

Now we follow Section 7.3 of \cite{GH}. Namely, 
let $\chi >0$ be sufficiently small that 
$\chi^\prime=e^{\chi LB}-1<\sqrt{2}-1$.
Note that $\chi $ only depends on $B$. 
For this number $\chi$ and for $x\in X$, 
$\alpha_1,\alpha_2\in {\cal A}(B)\times
{\cal A}(B)$ define
\begin{equation}\label{ultradist}
\tilde\delta_x(\alpha_1,\alpha_2)=
e^{-\chi\tau_B(x,\alpha_1,\alpha_2)}.\end{equation}
From Lemma \ref{tauestimate} and  
Proposition 7.3.10 of \cite{GH} we obtain.

\begin{corollary}\label{deltadistance}
There is a family $\{\delta_x\}$ $(x\in X)$ of distances on 
${\cal A}(B)$ with the following properties.
\begin{enumerate}
\item The topology on ${\cal A}(B)$ defined by the 
distances $\delta_x$ is the restriction of the
product topology on $\partial X\times\partial X-\Delta$.
In particular, $({\cal A}(B),\delta_x)$ is locally compact.
\item The distances $\delta_x$ are invariant under the involution
$\iota:(\xi,\eta)\to (\eta,\xi)$ of ${\cal A}(B)$ exchanging 
the two components 
of a point in ${\cal A}(B)$.
\item 
$(1-2\chi^\prime)\tilde \delta_x\leq \delta_x\leq \tilde \delta_x$
for all $x\in X$.
\item  
$e^{-\chi d(x,y)}\leq \delta_y\leq e^{\chi d(x,y)}\delta_x$ for all
$x,y\in X$.
\item The family $\{\delta_x\}$ is invariant under the action of 
${\rm Iso}(X)$ on ${\cal A}(B)\times X$. 
\end{enumerate}
\end{corollary}
\begin{proof} The existence of a family $\{\delta_x\}$
of distance functions with the property stated 
in the third part of the corollary is 
immediate from Lemma \ref{tauestimate} and Proposition 3.7.10 of \cite{GH}.
The second part of the corollary is immediate by invariance
of $\tau_B(x,\cdot,\cdot)$ under the involution $\iota$. The forth part 
follows from the 
construction of the distance $\delta_x$ from the functions
$\tilde \delta_x$ 
and from the third part of Lemma \ref{taubfirst}. Invariance
under the action of the isometry group is an immediate
consequence of invariance of the function $\tau_B$.

We are left with showing that for a given $x\in X$ the distance
$\delta_x$ induces the restriction of the product topology.
For this we first show that the identity
$({\cal A}(B),\delta_x)\to {\cal A}(B)\subset
\partial X\times \partial X-\Delta$ is continuous. For this
let $\gamma:\mathbb{R}\to X$ be a $B$-contracting
geodesic connecting $\xi$ to $\eta$.
By the definition of the distances $\delta_x$, 
if $(\xi_i,\eta_i)\to (\xi,\eta)$ in 
$({\cal A}(B),\delta_x)$ then there are $B$-contracting
geodesics $\gamma_i$ connecting $\xi_i$ to $\eta_i$ which
have longer and longer subsegments contained in a tubular
neighborhood of radius $6B+2$ about the geodesic 
$\gamma$. Moreover, the midpoints of these segments
are contained 
in a fixed compact subset of $X$. By the definition of the topology on 
$\partial X$, this implies that $(\xi_i,\eta_i)\to (\xi,\eta)$ in
$\partial X\times \partial X-\Delta$.

Continuity of the identity ${\cal A}(B)\subset \partial X\times \partial X
-\Delta
\to ({\cal A}(B),\delta_x)$ follows
in the same way. Namely, by the first part of the proof of 
Lemma \ref{abclosed}, if $(\xi_i,\eta_i)\subset {\cal A}(B)$,
if $(\xi_i,\eta_i)\to (\xi,\eta)\in \partial X\times \partial X
-\Delta$ with respect to the product topology and if
$\gamma_i$ is a $B$-contracting geodesic connecting
$\xi_i$ to $\eta_i$ then up to passing to a subsequence,
we may assume that the geodesics $\gamma_i$ converge
uniformly on compact sets to a $B$-contracting geodesic
$\gamma$ connecting $\xi$ to $\eta$. 
By convexity, by Lemma \ref{thintriangle} and by 
the definition of the function $\tau_B$, this implies  
that $(\xi_i,\eta_i)\to (\xi,\eta)$ in $({\cal A}(B),\delta_x)$
for every $x\in X$.
\end{proof}

The following analog
of Lemma 2.1 of \cite{H08b} is immediate from
Corollary \ref{deltadistance}.

\begin{lemma}\label{productmetric}
${\cal A}(B)\times X$ admits a natural ${\rm Iso}(X)$-invariant
$\iota$-invariant 
distance function $\tilde d$ inducing the
product topology. There is a number $c>0$
such that for every $x\in X$, the restriction of $\tilde d$ to
${\cal A}(B)\times \{x\}$ satisfies
\[c\delta_x(\alpha,\beta)\leq \tilde d((\alpha,x),(\beta,x))\leq
\delta_x(\alpha,\beta)\forall \alpha,\beta\in {\cal A}(B).\]
\end{lemma}
\begin{proof} The proof is identical to the proof of
Lemma 2.1 of \cite{H08b}. \end{proof}

\section{Continuous bounded cocycles}

This section is devoted to the proof of the 
main technical result of this note.
For its formulation, let again $X$ be a proper
${\rm CAT}(0)$-space and let $G$ be a \emph{closed} 
non-elementary subgroup of the
isometry group of $X$ with limit set $\Lambda$. Then
$G$ is a locally compact $\sigma$-compact topological group.
Assume that $G$  
contains a rank one element. Let 
$T\subset \Lambda^3$ be the space of triples of pairwise distinct
points in $\Lambda$. By Lemma \ref{northsouth}, 
$T$ is a locally compact 
uncountable topological $G$-space without
isolated points. As in Section 3,
for a number $B>0$ denote by ${\cal A}(B)\subset\partial X\times
\partial X-\Delta$ the set of pairs of distinct points in 
$\partial X$ which can be connected by a $B$-contracting geodesic.
Let moreover
$T(B)\subset T$ be the set of triples $(a_1,a_2,a_3)\in T$ with the
additional property that $(a_i,a_{i+1})\in {\cal A}(B)$
$(1\leq i\leq 3$ and 
where indices are taken modulo three). 
By Lemma \ref{abclosed}, $T(B)$ is closed
subset of $T$ which is invariant under the diagonal action of $G$.

For a Banach-module $E$ for $G$  
define an \emph{$E$-valued continuous
bounded two-cocycle} for the action of $G$ on $T(B)$ to be a
continuous \emph{bounded} $G$-equivariant map $\omega:T(B)\to E$
which satisfies the following two properties.
\begin{enumerate}
\item 
For every permutation $\sigma$ of
the three variables, the \emph{anti-symmetry condition}
$\omega\circ
\sigma={\rm sgn}(\sigma)\omega$ holds.
\item For every quadruple $(a_1,a_2,a_3,a_4)$
of distinct points in $\Lambda$ such that 
$(a_i,a_j)\in {\cal A}(B)$ for $i\not=j$ 
the \emph{cocycle equality}
\begin{equation}\label{coc}
\omega(a_2,a_3,a_4)-\omega(a_1,a_3,a_4)+
\omega(a_1,a_2,a_4)-\omega(a_1,a_2,a_3)=0
\end{equation}
is satisfied.
\end{enumerate}

Every locally compact $\sigma$-compact topological group $G$ 
admits a left invariant locally finite Haar measure $\mu$. 
For $p\in (1,\infty)$ denote by $L^p(G\times G,\mu\times \mu)$ 
the Banach space of all functions on $G\times G$ which
are $p$-integrable with respect to the product measure $\mu\times\mu$.
The group $G$ acts continuously and isometrically on 
$L^p(G\times G,\mu\times \mu)$ by left translation via
$(gf)(h,u)=f(gh,gu)$.

Let $C_b(G\times G)$ be the space of continuous bounded
functions on $G\times G$ and let $x_0\in X$ be an arbitrary
fixed point.
The following result is a modification of Theorem 2.3 of \cite{H08b}.

\begin{theo}\label{cocycle}
Let $X$ be a proper ${\rm CAT}(0)$-space 
and let $G< {\rm Iso}(X)$ be a non-elementary closed
subgroup.
For $B>0$, for every $p\in (1,\infty)$ and 
for every triple $(a,b,\xi)\in
T(B)$ such that $(a,b)$ is
the pair of fixed points of a rank-one element of $G$
there is a continuous map
$\alpha:{\cal A}(B)\to C_b(G\times G)$ with the following 
properties.
\begin{enumerate}
\item $g\circ \alpha(g\xi,g\eta)
=\alpha(\xi,\eta)=-\alpha(\eta,\xi)$ for all 
$(\xi,\eta)\in {\cal A}(B)$ and all $g\in G$.
\item For every $(\xi,\eta)\in {\cal A}(B)$ and all neighborhoods
$A_1$ of $\xi$, $A_2$ of $\eta$ in $X\cup \partial X$ the 
intersection of the support of $\alpha(\xi,\eta)$ with
the set $\{(g,h)\in G\times G\mid gx_0\in X-(A_1\cup A_2)\}$ is compact.
\item For every $p\in (1,\infty)$ 
the assignment \[\omega:(\sigma,\eta,\beta)\in T(B)\to 
\omega(\sigma,\eta,\beta)=
\alpha(\sigma,\eta)+\alpha(\eta,\beta)+\alpha(\beta,\sigma)\] 
is an
$L^p(G\times G,\mu\times \mu)$-valued continuous bounded
two-cocycle for the action of $G$ on $T(B)$.
\item $\omega(a,b,\xi)\not=0$.
\item If $G$ does not act transitively on $\Lambda\times\Lambda-\Delta$
and if $(a_i,b_i)\in {\cal A}(B)$ $(i=1,\dots,k)$ 
are pairs of fixed points of rank-one elements of
$G$ such that
the $G$-orbits of $(a_i,b_i),(a,b)$ are disjoint 
then we can choose 
$\alpha$ in such a way that the support of $\omega$
does not contain
$(a_i,b_i,\sigma)$ for $i=1,\dots,k$ and all $\sigma$.
\end{enumerate}
\end{theo}
\begin{proof}
Let $G< {\rm Iso}(X)$ be a closed non-elementary
subgroup which
contains a $B$-rank-one element for some $B>0$. 
We divide the proof of the theorem into five steps.

{\sl Step 1:}

Let $x_0\in X$ be an arbitrary point
and denote by $G_{x_0}$ the stabilizer of $x_0$ in $G$. Then
$G_{x_0}$ is a compact subgroup of $G$, and the quotient space
$G/G_{x_0}$ is $G$-equivariantly 
homeomorphic to the orbit $Gx_0\subset X$ of $x_0$.
Note that $Gx_0$ is a closed subset of $X$ and hence it is locally compact.
The group $G$ acts on the locally compact space
${\cal A}(B)\times Gx_0$ as a group of homeomorphisms.

The ${\rm Iso}(X)$-invariant metric $\tilde d$ on 
${\cal A}(B)\times X$ 
constructed in Lemma \ref{productmetric} induces a
$G$-invariant metric on ${\cal A}(B)\times G/G_{x_0}$ 
which defines the product topology. Hence we
obtain a 
$G$-invariant symmetrized product metric $\hat d$ on 
\begin{equation}V={\cal A}(B)\times G/G_{x_0}\times G/G_{x_0} \notag
\end{equation} 
by defining  
\begin{equation}\label{productmetric2}
\hat d((\xi,x,y), (\xi^\prime,x^\prime,y^\prime))
=\frac{1}{2}\bigl(\tilde d((\xi,x),(\xi^\prime,x^\prime))+
\tilde d((\xi,y),(\xi^\prime,y^\prime))).
\end{equation}
The topology defined on $V$ by this metric is 
the product topology, in particular it is locally compact.

Since $V$ is a locally compact $G$-space, 
the quotient space $W=G\backslash V$ admits a natural metric $d_0$ as
follows. Let 
\begin{equation}P:V\to W \notag
\end{equation}
be the canonical projection and define
\begin{equation}\label{d0}
d_0(x,y)=\inf\{\hat d(\tilde x,\tilde y)\mid P\tilde x=x, P\tilde
y=y\}. \end{equation}
The topology induced by this metric is the quotient
topology for the projection $P$. In particular, $W$ is a
locally compact metric space. A set
$U\subset W$ is open if and only
if $P^{-1}(U)\subset V$ is open. In other words,
open subsets of $W$ correspond precisely to 
$G$-invariant open subsets of $V$.
The projection $P$ is open and distance non-increasing.

The distance $\tilde d$ on ${\cal A}(B)\times X$ is invariant
under the involution 
$\iota:(\xi,\eta,x)\to (\eta,\xi,x)$ exchanging the two 
components of 
a point in ${\cal A}(B)$ and hence the same is true 
for the distance $\hat d$ on $V$.
Since the action of $G$ commutes
with the isometric involution $\iota$, the map  
$\iota$ descends to an isometric involution 
of the metric space $(W,d_0)$
which we denote again by $\iota$.

An open subset $U$ of $W$ is said to
have \emph{property} $(R_1,R_2)$ for some $R_1,R_2>0$ 
if for every 
$((\xi,\eta),gx_0,hx_0)\in P^{-1}(U)\subset V$ 
the distance in $X$ between $gx_0,hx_0$ and any geodesic
in $X$ connecting $\xi$ to $\eta$ is at most $R_1$ and
if moreover $d(gx_0,hx_0)\leq R_2$.

We claim that for every $w\in W$ there are numbers
$R_1,R_2>0$ and there is 
a neighborhood of $w$ in $W$ which has property $(R_1,R_2)$. 
Namely, let $v=((\xi,\eta),gx_0,hx_0)\in P^{-1}(w)$. 
Then $\xi$ can be connected to $\eta$ by a $B$-contracting 
geodesic $\gamma$ and therefore any geodesic
connecting $\xi$ to $\eta$ is contained in the 
$B$-tubular neighborhood of $\gamma$. 
By the discussion in Section 3 (see the proof of
Lemma \ref{abclosed}), there is a neighborhood $A$
of $(\xi,\eta)$ in ${\cal A}(B)$ such that
for all $(\xi^\prime,\eta^\prime)\in A$, any geodesic
connecting $\xi^\prime$ to $\eta^\prime$ passes through
a fixed compact neighborhood of $gx_0$. Thus  
by continuity, there are numbers
$R_1>0,R_2>0$ and there is an open 
neighborhood $U^\prime$ of $v$ in $V$ such that
for every $((\xi^\prime,\eta^\prime),g^\prime x_0,h^\prime x_0)\in U^\prime$
the distance between $g^\prime x_0,h^\prime x_0$ and any 
geodesic connecting $\xi^\prime$ to $\eta^\prime$ is at most $R_1$ 
and that moreover $d(g^\prime x_0,h^\prime x_0)\leq R_2$. 
However, distances and geodesics are preserved under isometries
and hence every point in $\tilde U=\cup_{g\in G}gU^\prime$
has this property. Since $\tilde U$ is open, $G$-invariant and 
contains $v$, the set $\tilde U$ 
projects to an open neighborhood of 
$w$ in $W$. This neighborhood has property $(R_1,R_2)$.

{\sl Step 2:} 

In equation (\ref{d0})
in Step 1 above, we defined a distance $d_0$ on the space
$W=G\backslash V$.
With respect to this distance, 
the involution $\iota$ acts non-trivially
and isometrically. Choose a small closed metric 
ball $D$ in $W$  
which is disjoint from its image under $\iota$. 
In Step 5 below we 
will construct explicitly 
such balls $D$, however 
for the moment, we simply assume that such a ball exists.
By Step 1 above, we may assume that $D$ has property $(R_1,R_2)$
for some $R_1,R_2>0$.

Let ${\cal H}$ be
the vector space of all H\"older continuous functions $f:W\to
\mathbb{R}$ supported in $D$. 
An example of such a function can be obtained as follows.

Let $z$ be an interior point of $D$ and let $r>0$ be sufficiently
small that the closed metric 
ball $B(z,r)$ of radius $r$ about $z$ is contained
in $D$. Choose a smooth function $\chi:\mathbb{R}\to [0,1]$ such that
$\chi(t)=1$ for $t\in (-\infty,r/2]$ and $\chi(t)=0$ for 
$t\in [r,-\infty)$ and 
define $f(y)=\chi(d_0(z,y))$. Since the function $y\to d_0(z,y)$ 
on $W$ is one-Lipschitz and $\chi$ is smooth, the function 
$f$ on $W$ is Lipschitz, does not vanish at $z$ 
and is supported in $D$.  

Since $D$ is disjoint from $\iota(D)$ by assumption and since
$\iota$ is an isometry, every function $f\in {\cal H}$ admits
a natural extension to a H\"older 
continuous function $f_0$ 
on $W$ supported in $D\cup \iota(D)$ whose restriction to
$D$ coincides with the restriction of $f$ and 
which satisfies $f_0(\iota z)=-f_0(z)$ for
all $z\in W$. 
The function $\hat f=f_0\circ P:V\to \mathbb{R}$
is invariant under the action of $G$, and it 
is \emph{anti-invariant} under the involution 
$\iota$ of $V$, i.e. it satisfies
$\hat f(\iota(v))=-\hat f(v)$ for all $v\in V$ (here as before,
$P:V\to W$ denotes the canonical projection).

Equip $\tilde V={\cal A}(B)\times G\times G$
with the product topology. The group $G$
acts on $G\times G$ by left translation, and it acts 
diagonally on $\tilde V$.
There is a natural continuous 
projection $\Pi:\tilde V\to V$ which is
equivariant with respect to the action of $G$ and with respect to
the action of the involution $\iota$ on $\tilde V$ and $V$.
The function $\hat f$ on $V$ lifts to
a $G$-invariant $\iota$-anti-invariant
continuous function $\tilde f=\hat f\circ \Pi$ on
$\tilde V$.

For $(\xi, \eta)\in {\cal A}(B)$ write 
\[F(\xi,\eta)=\{(\xi,\eta,z)\mid
z\in G\times G\}.\] 
The sets $F(\xi,\eta)$ define a $G$-invariant
foliation ${\cal F}$ of $\tilde V$. The leaf 
$F(\xi,\eta)$ of ${\cal F}$
can naturally be identified with $G\times G$. 
For all $(\xi,\eta)\in
{\cal A}(B)$ and every function $f\in {\cal H}$ we denote
by $f_{\xi,\eta}$ 
the restriction of the function $\tilde f$ to $F(\xi,\eta)$,
viewed as a continuous bounded function on
$G\times G$. 
For every $f\in {\cal H}$, all $(\xi,\eta)\in {\cal A}(B)$ 
and all $g\in G$
we then have $f_{g\xi,g\eta}\circ g=f_{\xi,\eta}=-f_{\eta,\xi}$.
Since the functions $f_{\xi,\eta}$ are restrictions 
to the leaves of the foliation ${\cal F}$ of a globally continuous
bounded function on $\tilde V$, 
the assignment 
\[(\xi,\eta)\in {\cal A}(B)\to \alpha(\xi,\eta)=
f_{\xi,\eta}\in C_b(G\times G)\]
is a continuous map of ${\cal A}(B)$ into  
$C_b(G\times G)$. By construction, it  
satisfies 
\[g\circ \alpha(g\xi,g\eta)=\alpha(\xi,\eta)=-\alpha(\eta,\xi)\,\forall
(\xi,\eta)\in {\cal A}(B),\forall g\in G.\] 
Thus the map $\alpha$ fulfills the 
first requirement in the
statement of the theorem.
The second requirement is also satisfied since 
the set $D$ is assumed to have property $(R_1,R_2)$
and since moreover for every geodesic $\zeta:\mathbb{R}\to X$, 
for every open neighborhood $A$ of $\zeta(\infty)\cup\zeta(-\infty)$ 
in $X\cup \partial X$ and for every $R>0$ 
the intersection of the closed $R$-neighborhood of $\zeta$ with 
$X-A$ is compact.

{\sl Step 3:}

For $f\in {\cal H}$ and for an ordered
triple $(\xi,\eta,\beta)\in T(B)$
define
\begin{equation}\label{omega}
\omega(\xi,\eta,\beta)= f_{\xi,\eta}+ 
f_{\eta,\beta}+
f_{\beta,\xi}\in C_b(G\times G).
\end{equation}
Since $f_{\xi,\eta}=-f_{\eta,\xi}$ 
for all $(\xi,\eta)\in {\cal A}(B)$, 
we have 
\[\omega\circ \sigma=({\rm
sgn}(\sigma))\omega\] for every permutation $\sigma$ of the three
variables. As a consequence, the cocycle condition 
for $\omega$ is also satisfied.
The assignment $(\xi,\eta,\beta)\in T(B)\to 
\omega(\xi,\eta,\beta)\in C_b(G\times G)$ is continuous with respect to the
compact open topology on $C_b(G\times G)$. Moreover, it is
equivariant with respect to the natural action of $G$ on 
the space $T(B)$ and on $C_b(G\times G)$.
This means that $\omega$ is a continuous bounded cocycle 
for the action of $G$ on $T(B)$ with
values in $C_b(G\times G)$.

For the proof of the theorem, we have to show 
that $\omega(\xi,\eta,\beta)\in L^p(G\times G,
\mu\times \mu)$ for every $p\in (1,\infty)$,
with $L^p$-norm bounded from above by a constant
which does not depend on $(\xi,\eta,\beta)$. For this 
let $(\xi,\eta,\beta)\in T(B)$ and let 
$\gamma:\mathbb{R}\to X$ be a 
$B$-contracting geodesic connecting $\xi$ to $\eta$.
By Lemma \ref{thintriangle2} there
is a point $y_0\in X$ which is contained in the
$\kappa_0=\kappa_0(B)=9B+6$-neighborhood 
of every side of a geodesic triangle with
vertices $\xi,\eta,\beta$ and side $\gamma$.
Via reparametrization of $\gamma$ we may assume that
$d(\gamma(0),y_0)\leq \kappa_0$.
Lemma \ref{thintriangle}, applied to the 
geodesic ray $\gamma[0,\infty)$ and a geodesic
$\zeta$ connecting $\beta$ to $\eta$ shows that there
is a subray of $\zeta$ whose Hausdorff distance to
$\gamma[2\kappa_0,\infty)$ is bounded from above
by $3B+1$.
Then by the
definition of the distances $\delta_x$ on ${\cal A}(B)$ and 
by Lemma \ref{decay} and Corollary \ref{deltadistance}, 
there is a number $r_0>0$ depending on $\kappa_0,R_1,R_2$ such that
if $t\geq 0$ and if $y\in X$
satisfies $d(\gamma(t),y)<\kappa_0+R_1+R_2$ then 
\begin{equation}\label{deltadifference}
\delta_y((\xi,\eta),(\beta,\eta)) \leq r_0e^{-\chi t}
\end{equation} 
where $\chi>0$ is as in Corollary \ref{deltadistance}.
Moreover, by the definition (\ref{productmetric2}) 
of the distance function $\hat d$ on $V$ and by
the estimate in Lemma \ref{productmetric} for the distance
function $\tilde d$ on ${\cal A}(B)\times X$, 
we have 
\begin{align}\label{hatdestimate}
\hat d\bigl((\xi,\eta,ux_0,hx_0),(\beta,\eta,ux_0,hx_0)\bigr)\\
\leq\frac{1}{2}(\delta_{ux_0}((\xi,\eta),(\beta,\eta))
+\delta_{hx_0}((\xi,\eta),(\beta,\eta)))
\leq r_0e^{-\chi t}\notag
\end{align}
whenever $d(ux_0,\gamma(t))\leq \kappa_0+R_1$ 
and $((\xi,\eta),ux_0,hx_0)$ is contained
in the support of $f_{\xi,\eta}$ or of $f_{\eta,\beta}$.

The function $\hat f=f_0\circ P:V\to \mathbb{R}$ 
constructed in Step 2 from the function $f\in {\cal H}$ 
is H\"older continuous and $\iota$-anti-invariant. Therefore
by the estimate (\ref{hatdestimate})
there are numbers $\sigma>0,r_1>r_0$ only depending on the H\"older
norm for $f$ with the following property.
Let $0\leq t$ and let $u,h\in G$ be such that
$d(ux_0,\gamma(t))<\kappa_0+R_1,d(hx_0,\gamma(t))<\kappa_0+R_1+R_2$; 
then\begin{equation}\label{hatf}
\vert \hat f(\xi,\eta,ux_0,hx_0)+
\hat f(\eta,\beta,ux_0,hx_0)\vert \leq r_1e^{-\sigma\chi t}.
\end{equation}
The function $f$ is bounded in absolute value
by a universal constant. Hence from 
the definition of the functions $f_{\xi,\eta}$ and
$f_{\eta,\beta}$ and from the estimate (\ref{hatf})
we obtain the existence of a constant
$r>r_1$ (depending on the H\"older norm of $f$) 
such that
\begin{equation}\label{asymptotic}
\vert (f_{\xi,\eta}+
f_{\eta,\beta})(u,h)\vert \leq re^{-\sigma\chi t}
\end{equation}
whenever $d(ux_0,\gamma(t))\leq \kappa_0+R_1$.

{\sl Step 4:}

Let $\nu=\mu\times \mu$ be the left invariant product measure on
$G\times G$. Our goal is to show that for every $p\in (1,\infty)$, 
the cocycle
$\omega$ defined in equation (\ref{omega}) above 
is in fact a bounded cocycle with values in
$L^p(G\times G,\nu)$. 
For this 
we show that for every $(\xi,\eta,\beta)\in T(B)$ and for every $p>1$
the $L^p$-norm of 
the function $\omega(\xi,\eta,\beta)$ with respect to the measure
$\nu$ on $G\times G$ is uniformly bounded and that
moreover the assignment 
$(\xi,\eta,\beta)\to \omega(\xi,\eta,\beta)\in L^p(G\times G,\nu)$ is
continuous.

For a subset $C$ of $X$ write 
\[C_{G,R_2}=\{(u,h)\in G\times G\mid
ux_0\in C,d(ux_0,hx_0)\leq R_2\}.\] 
We claim that 
there is a number $m>0$ such that for every subset
$C$ of $X$ of diameter at most $2R_1+4\kappa_0+1$ the 
$\nu$-mass of the set $N(C)_{G,R_2}$ 
is at most $m$. 
Namely, the set 
\[D=\{(u,h)\in G\times G\mid d(ux_0,x_0)\leq 2R_1+4\kappa_0+1,
d(ux_0,hx_0)\leq R_2\}\] 
of $G\times G$ is compact and hence
its $\nu$-mass is finite, say this mass equals $m>0$. 
On the other hand, if $C\subset X$ is a set of
diameter at most 
$2R_1+4\kappa_0+1$ and if there is some 
$g\in G$ such that $gx_0\in C$ then
any pair $(u,h)\in N(C)_{G,R_2}$ is contained
in $gD$. Our claim now
follows from the fact that $\nu$ is invariant under left
translation.

As in Step 3 above, let $(\xi,\eta,\beta)\in T(B)$ and
let $T$ be an ideal geodesic triangle with $B$-contracting
sides and vertices $\xi,\eta,\beta$. Let 
$y_0\in X$ be a point which is contained in the
$\kappa_0$-neighborhood of every side of $T$.
Let $\gamma:\mathbb{R}\to X$ be the side of $T$ connecting
$\xi$ to $\eta$, parametrized in such a
way that $d(y_0,\gamma(0))\leq \kappa_0$.
Also, let $\rho:\mathbb{R}\to X$ 
be the side of $T$ connecting $\beta$ to $\eta$ which is
parametrized in such a way that $d(y_0,\rho(0))\leq \kappa_0$. 
Then $\gamma[0,\infty),\rho[0,\infty)$ are two sides of 
a geodesic triangle in $X$ with vertices $\gamma(0),\rho(0),\eta$.
Since $d(\gamma(0),\rho(0))\leq 2\kappa_0$, by 
convexity of the distance function we have 
$d(\gamma(t),\rho(t))\leq 2\kappa_0$ for all $t\geq 0$.
In particular, the $R_1$-neighborhood of $\rho[0,\infty)$
is contained in the $R_1+2\kappa_0$-neighborhood of
$\gamma[0,\infty)$.

By construction, if $N(C,r)$ denotes the
$r$-neighborhood of a set $C$ then the 
support of the function $f_{\xi,\eta}$ is contained
in $N(\gamma(\mathbb{R}),R_1)_{G,R_2}$ and similarly
for the functions $f_{\eta,\beta},
f_{\beta,\xi}$.
As a consequence, 
the support of the function $\omega$ defined in
(\ref{omega}) above is the disjoint union of the 
three sets 
\begin{align}N(\gamma[0,\infty),R_1+2\kappa_0)_{G,R_2},&
N(\gamma(-\infty,-0],R_1+2\kappa_0)_{G,R_2},\\
N(\rho(-\infty,-0]&  ,R_1+2\kappa_0)_{G,R_2}.\notag\end{align}
Moreover, there is a number $\tau>0$ only depending on $R_1,R_2$
and $\kappa_0$ such that
the restriction of $\omega$ to
$N(\gamma[\tau,\infty),R_1+2\kappa_0)_{G,R_2}$ coincides with
the restriction of the function
$f_{\xi,\eta}+f_{\eta,\beta}$
and similarly for the other two sets in the
above decomposition of the support of $\omega$. 
Since $\omega$ is uniformly bounded, to show
that $\omega$ is contained in $L^p(G\times G,\nu)$ it is now enough to
show that there is constant $c_p>0$ only depending on $p$ and 
the H\"older norm of $f$ such that
\begin{equation}\label{integralestimate1}
\int_{N(\gamma[\tau,\infty),R_1+2\kappa_0)_{G,R_2}}\vert
f_{\xi,\eta}+ 
f_{\eta,\beta}\vert^p d\nu < c_p. \notag
\end{equation}

However, this is immediate from the 
estimate (\ref{asymptotic}) 
together with the control on
the $\nu$-mass of subsets of   
$N(\gamma[\tau,\infty),R_1+2\kappa_0)_{G,R_2}$.
Namely, we showed that for every integer $k\geq 0$ the
$\nu$-mass of the set
$N(\gamma[\tau+k,\tau+k+1],R_1+2\kappa_0)_{G,R_2}$ is bounded from
above by a universal constant $m>0$. 
Moreover, for every $p\geq 1$ the value of the function
$\vert f_{\xi,\eta}+f_{\beta,\eta}\vert^p $
on this set does not exceed $r^pe^{-p\sigma\chi(\tau+k)}$.
Thus the inequality \label{integralestimate} holds true with
$c_p=r^p\sum_{k=0}^\infty e^{-p\sigma\chi(\sigma+k)}$.

Since the function $\tilde f$ on $\tilde V$ is globally continuous,
the same consideration also
shows that $\omega(\xi,\eta,\beta)\in L^p(G\times G,\nu)$
depends continuously on $(\xi,\eta,\beta)\in T(B)$. 
Namely, let $(\xi_i,\zeta_i,\eta_i)\subset T(B)$ be a sequence
of triples of pairwise distinct points converging to a 
triple $(\xi,\eta,\beta)\in T(B)$. By the above consideration,
for every $\epsilon >0$ there is a compact subset $A$ of $G\times G$
such that $\int_{G\times G-A}\vert \omega(\xi_i,\eta_i,\beta_i)\vert^p
d\nu\leq \epsilon$ for all sufficiently large $i>0$ and that
the same holds true for $\omega(\xi,\eta,\zeta)$. 
Let $\chi_A$ be the characteristic function of $A$.
By continuity of the
function $\tilde f$ on $\tilde V$ and compactness,
the functions $\chi_A\omega(\xi_i,\eta_i,\beta_i)$ converge
as $i\to \infty$ in $L^p(G\times G,\nu)$ to $\chi_A\omega(\xi,\eta,\zeta)$.
Since $\epsilon >0$ was arbitrary, the required continuity follows.

By construction, the assignment
$(\xi,\eta,\beta)\to \omega(\xi,\eta,\beta)$ 
is equivariant under the action of $G$
on the space $T(B)$
and on $L^p(G\times G,\nu)$ and satisfies the
cocycle equality (\ref{coc}).
In other words, $\omega$ defines a continuous 
$L^p(G\times G,\nu)$-valued bounded cocycle for the action of $G$
on $T(B)$ as required.

{\sl Step 5:}

Let $g\in G$ be a $B$-rank one isometry, let $a\not=b\in \partial X$
be the attracting and repelling 
fixed point for the action of $g$ on $\partial X$, respectively,
and let $\xi\in \Lambda-\{a,b\}$ be such that
$(a,b,\xi)\in T(B)$. We have to show
that we can find a cocycle $\omega$ as in (\ref{omega}) above with
$\omega(a,b,\xi)\not=0$.

For this let $\gamma$ be a $B$-contracting oriented axis for $g$ and 
let $x_0\in \pi_{\gamma(\mathbb{R})}(\xi)$ be
the basepoint
for the above construction. 
The orbit of $x_0$ under the infinite cyclic subgroup of 
$G$ generated by $g$ is contained in
the geodesic $\gamma$.
Since $g^jx_0\to a,g^{-j}x_0\to b$ $(j\to \infty)$, there
are numbers 
$k<\ell$, $R_1 >2\kappa_0$ such that the $R_1+\kappa_0$-neighborhood
of the geodesic connecting $a$ to $\xi$ and the
$R_1+2\kappa_0$-neighborhood of the 
geodesic connecting $b$ to $\xi$ contains at most one of the
points $g^kx_0,g^\ell x_0$ and that moreover the distance
between $g^kx_0,g^\ell x_0$ is at least $4\kappa_0$. Choose 
$R_2>2d(g^kx_0,g^\ell x_0)$.

We claim that there is no $h\in G$ with 
$hg^kx_0=g^kx_0,hg^\ell x_0=g^\ell x_0$ and $h(a)=b,h(b)=a$. 
Namely, any isometry $h$ which exchanges $a$ and $b$ and fixes
a point on the axis $\gamma$ of $g$, say the point
$\gamma(t)$, maps the geodesic ray $\gamma[t,\infty)$ to
the geodesic ray $\gamma(-\infty,t]$. 
Thus a fixed point of $h$ on $\gamma$ is unique.
In particular, the 
projection of $(a,b,g^kx_0,g^\ell x_0)$ 
into the space $W$ is not fixed by the involution $\iota$.

The discussion at the end of Step 1 above shows that
we can find a function $f\in {\cal
H}$ supported in a ball $D\subset W$ about the projection 
of $(a,b,g^kx_0,g^\ell x_0)$ with property $(R_1,R_2)$ 
whose lift $\tilde f$ to $\tilde V$ 
does not vanish at
$(a,b,g^k,g^\ell)$. By the choice of $R_1$,
this means that $f_{a,b}(g^k,g^\ell)\not=0$
and $f_{b,\xi}(g^k,g^\ell)=
f_{\xi,a}(g^k,g^\ell)=0$. In other words, the
cocycle $\omega$ constructed as above
from $f$ does not vanish at $(a,b,\xi)$.

By Lemma \ref{closed} the $G$-orbit
of a pair of fixed points of a rank-one element of
$G$ is closed in $\Lambda\times \Lambda-\Delta$.
If the $G$-orbits of $(a_i,b_i),(a,b)$ 
$(i=1,\dots,k)$ are mutually disjoint where 
$(a_i,b_i)\in {\cal A}(B)$ are pairs of fixed points of such
rank-one elements then we can choose 
the function $f$ in such a way
that the support of its lift to $V$ does not intersect
the leaves of the foliation ${\cal F}$ determined by
the $G$-orbits of $(a_i,b_i)$. 
Thus we can find some cocycle $\omega$ which 
satisfies the fifth property of the theorem as well.
This completes the proof of the theorem.
\end{proof}

\section{Second continuous bounded cohomology}

Let $X$ be a proper ${\rm CAT}(0)$-space with isometry
group ${\rm Iso}(X)$. 
In this section we use Theorem \ref{cocycle} to construct 
nontrivial second bounded cohomology classes for 
closed non-elementary subgroups $G$ of ${\rm Iso}(X)$
with limit set $\Lambda$ 
which contain a rank-one element and 
act transitively on
the complement of the diagonal in $\Lambda\times \Lambda$.

We use the arguments from Section 3 of \cite{H08b}
(see also \cite{MMS04,MS04} for earlier results along
the same line).
Namely, let 
${\cal P}(\Lambda)$ be the space of all probability
measures on $\Lambda$, equipped with the
weak$^*$-topology. Denote moreover
by ${\cal P}_{\geq 3}(\Lambda)\subset {\cal P}(\Lambda)$ 
the set of all probability
measures which are not concentrated on at most two points.
We first show.

\begin{lemma}\label{tame}
Let $G<{\rm Iso}(X)$ be a non-elementary 
closed subgroup with limit
set $\Lambda$. If $G$ contains a rank-one element
and acts transitively on the complement of the
diagonal in $\Lambda\times \Lambda$ then 
the action of $G$ on ${\cal P}_{\geq 3}(\Lambda)$
is tame with compact point stabilizers.
\end{lemma}
\begin{proof}
Let $G<{\rm Iso}(X)$ be a closed non-elementary group 
with limit set $\Lambda$ which acts
transitively on the complement of the
diagonal $\Delta$ in $\Lambda\times \Lambda$ and
contains a rank one element $g\in G$. Then there is
a number $B>0$ and for every pair
$(\xi,\eta)\in \Lambda\times\Lambda-\Delta$ 
there is a $B$-contracting geodesic $\gamma:\mathbb{R}\to 
X$. This geodesic is the image of an axis of $g$ under an 
element of $G$.

Let $T\subset \Lambda^3$ be the space of 
triples of pairwise distinct points in $\Lambda$.
If $\gamma:\mathbb{R}\to X$ is a $B$-contracting geodesic then
every other geodesic connecting $\gamma(-\infty)$ to 
$\gamma(\infty)$ is contained in the $B$-tubular neighborhood of $\gamma$.
Thus by Lemma \ref{thintriangle2}, 
for every triple $(a,b,c)\in T$ there is a point
$x_0\in X$ whose distance to any of the sides
of any geodesic triangle in $X$ with vertices $a,b,c$
is at most $10B+6$. The set $K(a,b,c)$ of all
points with this property is clearly closed. 
Lemma \ref{thintriangle2} shows that it is moreover
of uniformly bounded diameter. In other words,
$K(a,b,c)$ is compact and hence it
has a unique center $\Phi(a,b,c)\in X$ 
where a center of a compact set $K\subset X$
is a point $x\in X$ such that the radius of the smallest
closed ball about $x$ containing $K$ is minimal (see p.10 in \cite{BGS85}). 

This construction defines a map 
$\Phi:T\to X$ which is equivariant with respect to the action of $G$.
Moreover, it is continuous. 
Namely, if $(a_i^1,a_i^2,a_i^3)\to (a^1,a^2,a^3)$ in $T$ 
then by the discussion in the proof of Lemma \ref{abclosed}
there is a compact neighborhood $A$ of 
$K(a^1,a^2,a^3)$ such that for all sufficiently
large $i$, every geodesic connecting a pair of points
$(a_i^j,a_i^{j+1})$ passes through $A$.
Since $X$ is proper by assumption, up to passing to a subsequence
we may assume that the compact sets
$K(a_i^1,a_i^2,a_i^3)$ of uniformly bounded diameter
converge in the Hausdorff 
topology for compact subsets of $X$ to a compact set 
$K$. On the other hand, up to passing
to a subsequence and reparametrization, a sequence of geodesics 
$\gamma_i^j$ connecting $a_i^j$ to $a_i^{j+1}$
converge
as $i\to \infty$ locally uniformly to a geodesic $\gamma^j$ connecting
$a^j$ to $a^{j+1}$. This implies that
$K\subset K(a^1,a^2,a^3)$.

However, by invariance under the action of $G$,
every geodesic connecting $a^j$ to $a^{j+1}$ is a limit
as $i\to \infty$ of a sequence of 
geodesics connecting $a^j_i$ to $a_i^{j+1}$
and therefore 
$K=K(a,b,c)$. Since the map which associates to a compact 
subset of $X$ its center
is continuous with respect to the Hausdorff topology on 
compact subsets of $X$, we conclude that the map $\Phi$ is in fact continuous.
Since the action of $G$ on $X$ is proper,
Lemma 3.4 of \cite{A96} then shows that the action of $G$ on $T$ is proper
as well.

Now the group $G<{\rm Iso}(X)$ is closed and hence the
stabilizer in $G$ of a point in $X$ is compact. 
By equivariance, the 
point stabilizer in $G$ of some 
$\mu\in {\cal P}_{\geq 3}(\Lambda)$ is compact as well.
This completes the proof of the lemma.
\end{proof}

The next easy consequence of a result of Adams and 
Ballmann \cite{AB98} 
will be important for the proof of Theorem 1.
For later reference, recall that the closure
of a normal subgroup of a topological group $G$ is normal,
and the closure of an amenable subgroup of $G$ is
amenable (Lemma 4.1.13 of \cite{Z}).

\begin{lemma}\label{amenrad}
Let $G<{\rm Iso}(X)$ be a closed non-elementary subgroup
which contains a rank-one element. Then 
a closed normal amenable subgroup $N$ of $G$ is compact, and $N$ fixes the
limit set of $G$ pointwise. 
\end{lemma}
\begin{proof}
Let $G<{\rm Iso}(X)$ be a closed non-elementary group which
contains a rank-one element and let $N\lhd G$ be a closed normal
amenable subgroup. Since $N$ is amenable, either
$N$ fixes a point $\xi\in \partial X$ or $N$ fixes
a flat $F\subset X$ \cite{AB98}. 

Assume first that $N$ fixes a point $\xi\in \partial X$.
Since $N$ is normal in $G$, for every $g\in G$ the point
$g\xi$ is a fixed point for $gNg^{-1}= N$.
On the other hand, by Lemma \ref{northsouth}, the closure
in $\partial X$ of every orbit for the action of $G$ 
contains the limit set $\Lambda$ of $G$, and the 
action of $G$ on
$\Lambda$ is minimal. Therefore by continuity,
$N$ fixes $\Lambda$ pointwise. 
Then Lemma \ref{tame} shows that $N$ is compact and hence
has a fixed point $x\in X$ by convexity.

If $N$ fixes a flat $F\subset X$ then 
we argue in the same way.
Namely, the image of $F$ under an isometry of $X$ is a flat.
Let $a\not=b\in \Lambda$ be the attracting and repelling
fixed points, respectively, of a rank one element $g$ of $G$. 
Then $a,b$ are visibility points in $\partial X$ and hence
there is no flat in $X$ whose
boundary in $\partial X$ contains one of the points $a,b$. 
Thus the boundary $\partial F\subset \partial X$  
of $F$ is contained in $\partial X-\{b\}$ and consequently
$g^k\partial F\to \{a\}$ 
$(k\to \infty)$. But by the argument
in the previous paragraph, $N$ fixes $g^k\partial F=
\partial (g^kF)$ and therefore
$N$ fixes $a$ by continuity. In other words,
$N$ fixes a point in $\partial X$. The first part of 
this proof then shows that indeed $N$ is compact.
This shows the lemma.
\end{proof}

{\bf Remark:} 
Caprace and Monod (Theorem 1.6 of 
\cite{CM08}, see also \cite{CF08}) found geometric conditions which
guarantee that 
an amenable normal subgroup of a non-elementary
group $G$ of isometries of $X$ vanishes.
This however need not be true 
under the above more general assumptions.
A simple example is a space of the form $X={\bf H}^2\times X_2$
where ${\bf H}^2$ is the hyperbolic plane and where 
$X_2$ is a compact ${\rm CAT}(0)$-space 
whose isometry group $H$ is 
non-trivial group. Then 
any axial isometry of ${\bf H}^2$ acts as a rank-one
isometry on $X$. The compact 
group $H$ is a normal subgroup of the isometry
group of $X$.

\bigskip
As in \cite{H08b}
we use Lemma \ref{tame} and Lemma \ref{amenrad} to show.

\begin{proposition}\label{boundedcoho1}
Let $G<{\rm Iso}(X)$ be a non-elementary
closed subgroup with limit set $\Lambda$. 
If $G$ contains a rank one element
and acts transitively on the complement of 
the diagonal in $\Lambda\times \Lambda$ then
$H_{cb}^2(G,L^p(G,\mu))\not=0$ for every
$p\in (1,\infty)$.
\end{proposition}
\begin{proof}
A \emph{strong boundary} for a locally compact
topological group $G$ is a standard
Borel space $(B,\nu)$ with a probability measure $\nu$ and 
a measure class
preserving amenable action of $G$ which is
\emph{doubly ergodic} (we refer to \cite{M} for
a detailed explanation of the significance of a strong
boundary). A strong boundary exists for
every locally compact topological group $G$ \cite{K03}. 

Let $G<{\rm Iso}(X)$ be a non-elementary closed
subgroup with limit set $\Lambda$. Assume that
$G$ acts transitively on the complement of the
diagonal in $\Lambda\times \Lambda$. 
Since the action of $G$ on its strong
boundary $(B,\nu)$ is amenable,
there is a $G$-equivariant measurable 
Furstenberg map $\phi:(B,\nu)\to {\cal P}(\Lambda)$
\cite{Z}.
By ergodicity of the action of $G$ on $(B,\nu)$,
either the set of all $x\in B$ with 
$\phi(x)\in {\cal P}_{\geq 3}(\Lambda)$ has full mass or
vanishing mass. 

Assume that this set has full mass. 
By Lemma \ref{tame}, the action of $G$ on 
${\cal P}_{\geq 3}(\Lambda)$ is tame, with compact
point stabilizers. Thus $\phi$ induces a $G$-invariant
map $(B,\nu)\to {\cal P}_{\geq 3}(\Lambda)/G$ which
is almost everywhere constant by ergodicity. 
Therefore by changing the map $\phi$ on a set
of measure zero, we can assume that $\phi$ is an 
equivariant map $(B,\nu)\to G/G_{\mu}$ 
where $G_{\mu}$ is the stabilizer of a point in
${\cal P}_{\geq 3}(\Lambda)$ and hence it is compact.
Since the action of $G$ on $(B,\nu)$ is amenable,
the group $G$ is amenable
(see p.108 of \cite{Z}). By assumption, 
$G$ does not fix globally a point in $\partial X$.
Then the group $G$ fixes a flat $F$ in $X$ \cite{AB98}.
The boundary $\partial F$ of $F$ is closed
and $G$-invariant and hence by Lemma \ref{northsouth} 
it contains the limit set
$\Lambda$ of $G$. But this means that there is a rank-one
element in $G$ with an axis contained in $F$ which is
impossible.

As a consequence, the image under $\phi$ of $\nu$-almost every $x\in B$
is a measure supported on at most two points.
By Lemma \ref{tame}, 
the action of $G$ on the space of triples of 
pairwise distinct points is proper and hence the
assumptions in Lemma 23 of \cite{MMS04} are satisfied.
We can then use Lemma 23 of \cite{MMS04} as in the
proof of Lemma 3.4 of \cite{MS04} to conclude that 
the image under $\phi$ of almost every $x\in B$ 
is supported in a single point. In other words,
$\phi$ is a $G$-equivariant 
Borel map of $(B,\nu)$ into $\Lambda$. Note that 
since the action of $G$ on $\Lambda$ is minimal, 
by equivariance the support of the measure class
$\phi_*(\nu)$ is all of $\Lambda$.

Let $\mu$ be a Haar measure of $G$.
By invariance under the action of $G$, 
there is some $B>0$ such that $(a,b)\in {\cal A}(B)$
for all $(a,b)\in \Lambda\times\Lambda-\Delta$.
Thus by Theorem \ref{cocycle},
for every $p\in (1,\infty)$ 
there is a nontrivial bounded continuous
$L^p(G\times G,\mu\times \mu)$-valued cocycle $\omega$ on the
space of triples of pairwise distinct points in $\Lambda$.
Then the
$L^p(G\times G,\mu\times \mu)$-valued
$\nu\times\nu\times \nu$-measurable bounded cocycle $\omega\circ \phi^3$
on $B\times B\times B$ is non-trivial on a set of positive
measure. Since $B$ is a strong boundary for $G$, this cocycle then
defines a \emph{non-trivial} class in $H_{cb}^2(G,L^p(G\times
G,\mu\times \mu))$ (see \cite{M}). On the other hand, the isometric
$G$-representation space $L^p(G\times G,\mu\times \mu)$ is a
direct integral of copies of the isometric $G$-representation space
$L^p(G,\mu)$ and therefore by Corollary 2.7 of \cite{MS04} and
Corollary 3.4 of \cite{MS06}, if $H_{cb}^2(G,L^p(G,\mu))=\{0\}$
then also $H_{cb}^2(G,L^p(G\times G,\mu\times \mu))=\{0\}$. This
shows the proposition. \end{proof}

\section{Second bounded cohomology}

In this section we investigate non-elementary closed 
subgroups of ${\rm Iso}(X)$
with limit set $\Lambda$ which contain a rank-one element
and which 
do not act transitively on the complement of the diagonal
in $\Lambda\times\Lambda$. Such a group $G$ is a locally compact
$\sigma$-compact group which admits a Haar measure $\mu$.
Our goal is to show that
for every $p\in (1,\infty)$ the second bounded cohomology
group $H_{b}^2(G,L^p(G,\mu))$ is infinite dimensional.

Unlike in Section 5,  
for this we can not use Theorem \ref{cocycle} directly
since the cocycle constructed in this theorem
may not be defined on the entire space of triples
of pairwise distinct points in $\Lambda$.
Instead we use the strategy from the proof of Theorem \ref{cocycle}
to construct
explicitly for every $p\in (1,\infty)$ bounded cocycles for
$G$ with values in $L^p(G,\mu)$ which define an infinite
dimensional subspace of $H_{b}^2(G,L^p(G,\mu))$.
Unfortunately, our construction does not yield continuous
bounded cocycles, so we only obtain information on the
group $H_b^2(G,L^p(G,\mu))$. (The construction of Bestvina
and Fujiwara \cite{BF08} does not yield continuous
bounded cocycles with real coefficients either. However,
these cocycles are integer-valued, and 
by an observation of Caprace and Fujiwara \cite{CF08},
such cocycles define in fact continuous bounded
cohomology classes.)

A \emph{twisted $L^p(G,\mu)$-valued quasi-morphism}
for a closed subgroup
$G$ of ${\rm Iso}(X)$ is a map $\psi:G\to L^p(G,\mu)$ such that
\[\sup_{g,h}\Vert \psi(g)+g\psi(h)-\psi(gh)\Vert_p<\infty\]
where $\Vert \,\Vert_p$ is the $L^p$-norm for functions on $G$.

Every unbounded twisted $L^p(G,\mu)$-valued quasi-morphism
for $G$
defines a second bounded cohomology class in $H_b^2(G,L^p(G,\mu))$ which 
vanishes if and only if there is a cocycle
$\rho:G\to L^p(G,\mu)$ (i.e. $\rho$ satisfies the
cocycle equation $\rho(g)+g\rho(h)-\rho(gh)=0$) such that
$\psi-\rho$ is bounded (compare the discussion in 
\cite{H08a}). We use twisted quasimorphisms to 
complete the proof of Theorem \ref{thm2} from the
introduction.

\begin{proposition}\label{boundedcoho2}
Let $G<{\rm Iso}(X)$
be a closed non-elementary subgroup with limit set $\Lambda$ which contains a
rank-one element. 
If $G$ does not act transitively
on the complement of the diagonal 
in $\Lambda\times \Lambda$ then for every
$p\in (1,\infty)$ the second bounded cohomology
group $H_{b}^2(G,L^p(G,\mu))$ is infinite dimensional.
\end{proposition}
\begin{proof}
Let $G<{\rm Iso}(X)$ be a closed subgroup with limit set
$\Lambda\subset \partial X$ which contains a rank-one element
and which does not act transitively on the complement of
the diagonal $\Delta$ in $\Lambda\times \Lambda$.

Let $g\in G$ be a rank-one element 
with attracting and repelling 
fixed points $a,b\in \Lambda$. By Lemma \ref{free} we may assume that
there is no $u\in G$ with $u(a,b)=(b,a)$.
Let $B_0>0$ be such that
every geodesic in $X$ connecting $a$ to $b$ is $B_0$-contracting.
Such a number exists by Lemma \ref{localray}.
By Lemma \ref{compact}, we can find some $h\in G$ such that
$hb\not=b$ and that the stabilizer ${\rm Stab}(b,hb)$ 
in $G$ of the pair of points
$(b,hb)$ is compact. 

By the consideration in the proof of Lemma \ref{compact}, 
there is a number
$B>B_0$ depending on $B_0$ and $h$ such that every
geodesic connecting $b$ to $hb$ is $B$-contracting.
Since ${\rm Stab}(b,hb)$ is compact,  
there is a such a geodesic $\gamma$ which is fixed pointwise by 
${\rm Stab}(b,hb)$.
Let $\gamma^{-1}$ be 
the geodesic obtained by reversing the orientation
of $\gamma$ (we note that for all of our constructions,
only the orientation of a geodesic but not an explicit
parametrization plays any role).
By invariance under isometries, for every $k\in \mathbb{Z}$,
every geodesic connecting $b$ to $g^khb$ is $B$-contracting.

For this number $B>0$ let $C=C(B)>0$ be as in Lemma \ref{localray}.
The $G$-orbit of $b$ consists of visibility points. Thus
for $u\not=v\in G$ with $ub\not=vb$ there is 
a (parametrized) geodesic 
$\xi$ connecting $ub$ to $vb$. The geodesic 
$vh^{-1}\gamma$ connects
$vh^{-1}b$ to $vb$ and hence 
by Lemma \ref{localray}, the geodesic $\xi$ passes 
through the $9B+6$-neighborhood of every point in  
$\pi_{vh^{-1}\gamma(\mathbb{R})}(ub)$ (this also holds true
if $ub=vh^{-1}b$). 
If $b(u,v,\xi)\in \mathbb{R}\cup \{-\infty\}$ 
is the smallest number such that 
$\xi(b(u,v,\xi))$ is 
contained in the $9B+6$-neighborhood of a point in 
$\pi_{vh^{-1}\gamma(\mathbb{R})}(ub)$
then the
geodesic ray $\xi[b(u,v,\xi),\infty)$ is 
$C$-contracting. This ray only depends on 
the geodesic line $\xi$, on the
coset $[v]$ of $v$ in $G/{\rm Stab}(b,hb)$ 
and on $ub\in \partial X$.
If $ub=b$ and if $v=g^kh$ then we have 
$b(u,v,\xi)=-\infty$. 

Define similarly a number 
$a(u,v,\xi)\in \mathbb{R}\cup \{\infty\}$
using the above procedure for the inverse of the geodesic $\xi$ and
the geodesic $uh^{-1}\gamma$.
The resulting ray $\xi(-\infty,a(u,v,\xi)]$ 
depends on $[u]\in G/{\rm Stab}(b,hb)$ and on $vb\in \partial X$
(where here and in the sequel, $[u]$ denotes the coset
of $u\in G$ in $G/{\rm Stab}(b,hb)$).

Let ${\cal G}(C)$ be the set of all (oriented)
geodesics $\eta:\mathbb{R}\to X$ such that there is an 
open connected relatively compact (perhaps empty) set  
$(a(\eta),b(\eta))\subset \mathbb{R}$  
with the property that $\eta(\mathbb{R}-(a(\eta),b(\eta)))$ is
$C$-contracting. For $\eta\in {\cal G}(C)$ we consider 
the subarc $\eta(a(\eta),b(\eta))$ to be 
part of the structure of $\eta$. 
Thus the same geodesic with two distinct subarcs removed
defines two distinct points in ${\cal G}(C)$.
Sometimes we write
$(\eta,(a(\eta),b(\eta)))\in {\cal G}(C)$ to describe
explicitly the data which make out a point in ${\cal G}(C)$.
The group $G$ naturally acts on ${\cal G}(C)$ from the left.

The above construction associates to any ordered pair of points
$([u],[v])\in 
G/{\rm Stab}(b,hb)\times G/{\rm Stab}(b,hb)$
with 
$ub\not=vb$ and every geodesic $\xi$ connecting 
$ub$ to $vb$ a (possibly empty) 
subarc $\xi(a(u,v,\xi),b(u,v,\xi))$ of $\xi$ in such a way
that $(\xi,(a(u,v,\xi),b(u,v,\xi)))\in {\cal G}(C)$.
The assignment
\[\Pi:([u],[v],\xi)\to \Pi([u],[v],\xi)= 
(\xi,(a(u,v,\xi),b(u,v,\xi)))\in {\cal G}(C)\]
satisfies the following properties.
\begin{enumerate}
\item Coarse independence on the choice of geodesics:
If $\xi,\zeta$ are two 
geodesics connecting $ub$ to $vb$ 
then the Hausdorff distance between the rays\\
$\xi[b(u,v,\xi),\infty),\zeta[b(u,v,\zeta),\infty)$ and
the between the rays\\ $\xi(-\infty,a(u,v,\xi)]$,
$\zeta(-\infty,a(u,v,\zeta)]$ is bounded
from above by $18B+12$.
\item Invariance under the action of $G$:
For $q\in G$ and a geodesic $\xi$ connecting
$ub$ to $vb$ we have 
$\Pi([qu],[qv],q\xi)=q\Pi([u],[v],\xi)$.
\end{enumerate}

Let $\tilde {\cal A}(b)$ be the union 
of the set $\{\Pi([u],[v],\xi)\mid u,v\in G, \xi\}\subset{\cal G}(C)$ 
with the $G$-translates of all oriented geodesics
connecting $a$ to $b$ or connecting $b$ to $a$ 
viewed as elements of ${\cal G}(C)$ with an empty
subarc.
Call two geodesics $\xi,\zeta\in \tilde{\cal A}(b)$ 
which are obtained in the above way from
the same pair of points $[u],[v]\in G/{\rm Stab}(b,hb)$ 
or which connect
the same pair of points in the $G$-orbit of $(a,b)$
or of $(b,a)$ 
\emph{equivalent}. This clearly defines an equivalence relation
on $\tilde {\cal A}(b)$.
Let $\hat{\cal A}(b)$ be the set of all
equivalence classes 
in $\tilde {\cal A}(b)$. Such an equivalence class 
either is determined 
by an ordered pair 
$([u],[v])\in G/{\rm Stab}(b,hb)\times G/{\rm Stab}(b,hb)$
and will then be denoted by $[u,v]$, or it is determined
by an ordered pair of points $(ua,ub)$ or $(ub,ua)$ 
in the $G$-orbit
of $(a,b)$ or $(b,a)$ and will be denoted by $[a,u]$ or $[u,a]$,
respectively.
The group $G$ naturally acts on $\hat {\cal A}(b)$
from the left.

Recall from Section 3 the definition of the function $\tau_C$
which associates to a point $y\in X$ and two 
(finite or infinite) geodesics 
$\zeta_1,\zeta_2$ with at most one common
endpoint in $\partial X$ a number
$\tau_C(x,\zeta_1,\zeta_2)\geq 0$. 
For $x\in X$ and for two geodesics
$\gamma_1,\gamma_2\in {\cal G}(C)$ define
\begin{align}
\tau_{C{\rm rel}}(x,\gamma_1, & \gamma_2)=\notag \\
\max\{\tau_C(x,\gamma_1[b(\gamma_1),\infty),\gamma_2[b(\gamma_2),\infty)), &
\tau_C(x,\gamma_1(-\infty,a(\gamma_1)],\gamma_2(-\infty,a(\gamma_2)]).\notag
\end{align}
Note that if the empty subarc is associated to each of 
the geodesics $\gamma_1,\gamma_2$
then $\tau_{C{\rm rel}}(x,\gamma_1,\gamma_2)=
\tau_C(\gamma_1,\gamma_2)$.

For $[u,v],[w,z]\in \hat{\cal A}(b)$ and $x\in X$ define
\[\tau_{C{\rm rel}}(x,[u,v],[w,z])=
\inf\tau_{C{\rm rel}}(x,\gamma_1,\gamma_2)\]
where the infimum is taken over all geodesics 
$\gamma_1,\gamma_2\in {\cal G}(C)$
in the equivalence class of $[u,v],[w,z]$.
By construction, we have 
\[\tau_{C{\rm rel}}(x,[u,v],[w,z])=\tau_{C{\rm rel}}(x,[v,u],[z,w])\]
for all $(u,v),(z,w)\in \hat{\cal A}(b)$.

Lemma \ref{tauestimate} shows that the function
$\tau_{C{\rm rel}}$ on $\hat{\cal A}(b)$ satisfies
the ultrametric inequality. Thus as in Section 3,
for each $x\in X$ we can use the function 
$\tau_{C{\rm rel}}(x,\cdot,\cdot)$ to define  
a pseudo-distance 
$\{\hat\delta_x^{C{\rm rel}}\}$ on $\hat {\cal A}(b)$.
These pseudo-distance functions are mutually bilipschitz
equivalent and hence they induce a family of distance 
functions $\{\delta_x^{C{\rm rel}}\}$ on the quotient
${\cal A}(b)$ of $\hat{\cal A}(b)$ 
obtained by identifying points with vanishing
pseudo-distance. Note that if two points in $\hat{\cal A}(b)$
are identified then they correspond to pairs
$([u],[v]),([u^\prime],[v^\prime])\in 
G/{\rm Stab}(b,hb)\times G/{\rm Stab}(b,hb)$ with $ub=u^\prime b,
vb=v^\prime b$. By abuse of notation, we denote the
projection of $[u,v]\in \hat{\cal A}(b)$ again by
$[u,v]$. The family $\{\delta_x^{C{\rm rel}}\}$
is invariant under the natural
action of $G$ on $X\times {\cal A}(b)$, and it is invariant
under the natural involution $\iota$ defined
by $\iota[u,v]=[v,u]$. In the sequel we always equip 
${\cal A}(b)$ with the topology induced by one
(and hence each) of these distance functions. 

Since to each geodesic in $X$ connecting $a$ to $b$ or
connecting $b$ to $g^khb$ we associated the empty
subarc, 
for each $x\in X$ the points $[e,g^kh]$ 
converge as $k\to \infty$
in $({\cal A}(b),\delta^{C{\rm rel}}_x)$ to 
$[e,a]$ (here as usual, $e$ is the unit in $G$).
In particular, the point 
$[e,a]\in {\cal A}(b)$ is not isolated for $\delta_x^{C{\rm rel}}$.
We use the distances $\delta_x^{C{\rm rel}}$ 
as in Lemma \ref{productmetric} 
to construct a $G$-invariant 
distance $\rho$ on ${\cal A}(b)\times X$
with the properties stated in Lemma \ref{productmetric}. 
Then $\rho$ induces the product topology on ${\cal A}(b)\times X$.

We now use the strategy from the 
proof of Theorem \ref{cocycle}.
Namely, let $x_0\in X$ be a point on an axis for the
rank-one element $g\in G$. Let $G_{x_0}$ be the stabilizer
of $x_0$ in $G$ and let 
$V(b)={\cal A}(b)\times G/G_{x_0}={\cal A}(b)\times Gx_0$.
The group $G$ acts on $V(b)$ as a group of isometries with 
respect to the restriction of the distance $\rho$.
Define $W=G\backslash V(b)$ and let $P:V(b)\to W$
be the canonical projection. The distance $\rho$ on $V(b)$
induces a distance $\hat \rho$ on $W$ by defining
$\hat\rho(x,y)=\inf\{\rho(\tilde x,\tilde y)\mid 
P\tilde x=x,P\tilde y=y\}$. Note that we have
$\hat \rho(x,y)>0$ for $x\not=y$ by 
the definition of the distance $\rho$ and 
the fact that the distances $\{\delta_x^{C{\rm rel}}\}$ 
depend uniformly Lipschitz continuously on $x\in X$.

The isometric involution $\iota$ of $({\cal A}(b)\times X,\rho)$
decends to an isometric involution on $W$ again denoted
by $\iota$. 
Since there is no $u\in G$ with $u(a,b)=(b,a)$, we can
find an open neighborhood $D$ of 
$w=P([e,a],x_0)\in V(b)$ 
which is disjoint from
its image under $\iota$. We choose $D$ to be contained in the
image under the projection $P$ of the set
${\cal A}(b)\times K$ where $K$ is the closed ball of radius
$1$ about $x_0$ in $Gx_0\subset X$.

Let $C_b(G{x_0})$ be the Banach space of continuous bounded
functions on $G{x_0}\subset X$. 
As in the proof of Theorem \ref{cocycle}, we use the induced
distance on $W$ to construct from a H\"older continuous function
$f$ supported in $D$ with $f(w)>0$ 
a $G$-invariant $\iota$-anti-invariant 
uniformly bounded continuous map 
$\tilde\alpha:{\cal A}(b)\to 
C_b(G{x_0})$ which lifts to a bounded continuous map 
$\alpha: {\cal A}(b)\to C_b(G)$
with the equivariance
properties as stated in this theorem.

By invariance and the definition of the distances
$\delta_x^{C{\rm rel}}$,
if $z\in G$ is such that $zx_0$
is contained in the support of
the function $\tilde\alpha[u,v]$ 
then there is a geodesic
$\xi\in {\cal G}(C)$ connecting
$ub$ to $vb$ so that 
$ux_0$ is contained in a tubular neighborhood 
of $\xi(-\infty,a(u,v,\xi)]\cup \xi[b(u,v,\xi),\infty)$ of uniformly bounded
radius.

Let $A$ be a small compact neighborhood of $b$ in $X\cup \partial X$
which does not contain the attracting fixed point $a$ of $g$.
For $u\in G$ with $ub\not=b$ define
a function $\Psi_\alpha(u):G\to \mathbb{R}$ by
\[\Psi_\alpha(u)(w)=\alpha[e,u](w)\]
if $wx_0\in X-(A\cup uA)$, and define 
$\Psi_\alpha(u)(w)=0$ otherwise.
By the above construction, for every $u\in G$ 
the function $\Psi_\alpha(u)$ is
measurable and supported in a compact subset
of $G$. Moreover, it is pointwise uniformly 
bounded independent of $u$. 
For every compact subset $K_0$ of $G$ there is a compact
subset $C$ of $G$ containing the support of each of the
functions $\Psi_{\alpha}(u)$ $(u\in K_0)$.
In particular, 
we have $\Psi_\alpha(u)\in L^p(G,\mu)$
for every $p>1$, and for every compact subset $K_0$ of $G$
the set $\{\Psi_{\alpha}(u)\mid u\in K_0\}\subset 
L^p(G,\mu)$ is bounded. If $ub=b$ then define $\Psi_\alpha(u)\equiv 0$. 

We claim that $\Psi_\alpha$ is unbounded. For this
note that as $k\to \infty$ we have 
$[e,g^kh]\to [e,a]$ in ${\cal A}(b)$ and 
$g^khA\to \{a\}$. In particular,
if $\xi$ is a $B$-contracting axis for $g$ containing
the point $x_0$ then $X-A-g^khA$ contains longer and
longer subsegments of $\xi$ which uniformly fellow-travel
the geodesic $g^k\gamma$ connecting $b$ to $g^khb$.
Now the function $\alpha[e,a]$ is invariant under the 
action of the rank-one element $g$ and its support contains the
point $x_0$. This implies that $\alpha[e,a]$ is \emph{not}
integrable. But then for $p>1$ the $L^p$-norm of the
functions $\Psi_\alpha(g^kh)$ tends to infinity as $k\to \infty$.

Define a function
$\omega:G^3\to L^p(G,\mu)$ by
\begin{align}\omega(u,uw,uh) 
= \omega(e,w,h) &= \Psi_\alpha(w)
+ w\Psi_\alpha(h)
- \Psi_\alpha(wh)\notag\\
=\hat\alpha[e,w]
+ w\hat\alpha[e,h]-\hat\alpha[e,wh]&  =
\hat\alpha[e,w]+ \hat\alpha[w,wh]+\hat\alpha[wh,e]\notag
\end{align}
where $\hat\alpha[u,v]$ is the restriction of 
$\alpha[u,v]$ to $\{z\mid zx_0\in X-(uA\cup vA)\}$.
Then $\omega$ is invariant under the diagonal action of $G$,
and we have $\omega\circ \sigma={\rm sgn}(\sigma)\omega$
for every permutation of the three variables.
Moreover, $\omega$ satisfies the cocycle identity. 
In other words, for every $p\in (1,\infty)$, $\omega$ 
is an $L^p(G,\mu)$-valued 2-cocycle for $G$.

We claim that the image of
$\omega$ is uniformly bounded. For this we argue as in the proof of
Theorem \ref{cocycle}. Namely, by assumption, if
$\Psi_\alpha(v)(w)\not=0$ 
then the point $wx_0$ is contained in
a uniformly bounded neighborhood of the geodesic 
$vh^{-1}\gamma$ 
and hence by the definition 
of the subarc $\xi(a(u,v,\xi),b(u,v,\xi))$ of a geodesic $\xi$ connecting
$b$ to $vb$ and
the estimate \ref{decay}, the $L^p$-norm of 
the restriction of $\alpha[e,v]+\alpha[v,z]$
to a tubular neighborhood of 
$\xi[b(u,v,\xi),\infty)$ is uniformly bounded.
By symmetry and the properties of the
support of the functions $\alpha[e,v]$, this implies as in 
the proof of Theorem \ref{cocycle} that
$\omega$ is uniformly bounded.

As a consequence, $\Psi_\alpha$ defines an element
in the kernel of the natural homomorphism 
$H_{b}^2(G,L^p(G,\mu))\to 
H^2(G,L^p(G,\mu))$ which vanishes
if and only if there is a map 
$\eta:G\to L^2(G,\mu)$
with $\eta(gh)=\eta(g)+g\eta(h)$ for all $g,h\in G$ and such that
$\Psi_\alpha-\eta$ is bounded (compare the discussion in \cite{H08b}).

Now recall from Lemma \ref{free} that 
if $G$ does not act transitively on its limit set then 
$G$ contains
a free subgroup $\Gamma$ with two generators consisting of rank one elements
which contains elements from infinitely many
conjugacy classes of $G$. 
If $g,h$ are the generators for $\Gamma$, if
$\Psi_\alpha(g)=\Psi_\alpha(h)=0$ and if $\Psi_\alpha$ defines the 
trivial class in $H_{b}^2(G,L^p(G,\mu))$ then   
$\Psi_\alpha\vert \Gamma$ is bounded. However, 
since the $G$-orbit of any pair of fixed points of rank-one elements
in $G$ is a closed subset of $\Lambda\times \Lambda-\Delta$, as in the
proof of the 
fifth property in Theorem \ref{cocycle} we conclude that there are infinitely
many linearly independent distinct such classes which pairwise
can not be obtained from each other by adding a bounded function.
This shows the proposition.
\end{proof}

{\bf Remark:} The construction in the proof of
Proposition \ref{boundedcoho2} does not yield continuous
bounded cohomology classes since in general the topology
on ${\cal A}(b)$ induced by the distance
functions $\delta_x^{C{\rm rel}}$ does not coincide with 
the restriction of the product topology.
However, the twisted quasi-morphisms constructed in the course
of the proof take on uniformly bounded values on every
compact subset of the group $G$.

\section{Structure of the isometry group}

In this section we use the results from Section 5 
and Section 6 to complete the
proof of Theorem 1 from the introduction.

\begin{proposition}\label{firstpart} 
Let $X$ be a proper ${\rm CAT}(0)$-space 
and let $G< {\rm Iso}(X)$ be a closed
subgroup which contains a rank-one element.
Then one of the following three
possibilities holds.
\begin{enumerate}
\item $G$ is elementary.
\item  $G$ contains an open subgroup $G^\prime$ of
finite index which is a compact extension of 
a simple Lie group of rank 1.
\item $G$ is a compact extension of a totally disconnected
group.
\end{enumerate}
\end{proposition}
\begin{proof}
Let $G$ be a closed subgroup of the isometry group
${\rm Iso}(X)$ of a proper ${\rm CAT}(0)$-space $X$.
Then $G$ is locally compact. Assume that $G$ 
is non-elementary and contains a rank-one element.
Then by Lemma \ref{amenrad}, the maximal normal amenable
subgroup $N$ of $G$ is compact,
and the quotient $L=G/N$ is a locally compact
$\sigma$-compact group. Moreover, $N$ acts
trivially on the limit set $\Lambda$ of $G$.

By the solution to Hilbert's fifth problem (see Theorem 11.3.4 in
\cite{M}), after possibly replacing $L$ by an open subgroup of
finite index (which we denote again by $L$ for simplicity), the
group $L$ splits as a direct product $L=H\times Q$ where
$H$ is a semisimple connected Lie group with finite center and
without compact factors and $Q$ is totally disconnected. If $H$
is trivial then $G$ is a compact extension of a totally
disconnected group.

Now assume that $H$ is nontrivial. Let $H_0<G$ and
$Q_0<G$ be the preimage
of $H,Q$ under the projection $G\to L$.
Then $H_0$ is not compact and
the limit set $\Lambda_0\subset \Lambda$ 
of $H_0< G$ is nontrivial. Since
$Q$ commutes with $H$ and the group $N$ 
acts trivially on $\Lambda$, the group $Q_0$
acts trivially on $\Lambda_0$ (this is discussed in the proof of 
Proposition 4.3 of \cite{H08b}, 
and the proof given there is valid in our situation 
as well). In particular, if $\Lambda_0$ consists of a single point
then $G$ is elementary.
As a consequence, if $G$ is non-elementary then 
$Q_0$ fixes at least
two points in $\partial X$.

We show that $\Lambda_0=\Lambda$. Since $\Lambda_0\subset
\Lambda$ is closed, by Lemma \ref{closed} it suffices
to show that the fixed points of every rank-one element
of $G$ are contained in $\Lambda_0$. Thus
let $g\in G$ be a rank-one element. By Lemma \ref{northsouth},
$g$ acts with north-south dynamics on $\partial X$ with
attracting fixed point $a\in \Lambda$ and repelling
fixed point $b\in \Lambda$. 
If $a\not\in \Lambda_0$ then there is a point 
$\xi\in \Lambda_0-\{a,b\}$.
Write $g=g_0q$ with $g_0\in H_0,q\in Q_0$.
Since $g_0$ and $q$ commute up to a compact normal subgroup which
fixes $\Lambda\supset \Lambda_0$ pointwise, we have
$g_0^k\xi=g_0^kq^k\xi=g^k\xi\to a$ $(k\to \infty)$. 
But $g_0^k\xi\in \Lambda_0$ for all $k>0$
and therefore by compactness we have $a\in \Lambda_0$. 
Since $a$ was an arbitrary fixed point of a 
rank-one element in $G$ we conclude that $\Lambda_0=\Lambda$
and hence
$Q_0$ fixes the limit set of $G$ pointwise.
However, since $G$ is non-elementary by assumption, 
in this case the argument in 
the proof of Lemma \ref{amenrad} shows that
$Q_0$ is compact and hence $Q$ is trivial.

To summarize, if $G$ is non-elementary 
then up to passing to an
open subgroup of finite index, either $G$ is a compact extension of a  
totally disconnected group
or $G$ is a compact extension of a 
semisimple Lie group $H$ with finite center and without compact
factors.

We are left with showing that if the group 
$G$ is a compact extension of a
semisimple Lie group $H$ with finite center and without compact
factors then $H$ is simple and of rank 1.
For this assume first that $G$ acts transitively
on the complement of the
diagonal in $\Lambda\times \Lambda$.
Then Proposition \ref{boundedcoho1} shows that 
$H_{cb}^2(G,L^2(G,\mu))\not=\{0\}.$ 
By Corollary 8.5.2 of \cite{M}, this implies that
$H_{cb}^2(H,L^2(H,\mu))\not=\{0\}$ as well. 
By the super-rigidity result for bounded cohomology
of Burger and Monod \cite{BM99}, we conclude
that $H$ is simple of rank one (see \cite{MS04,H08b}
for details on this argument).

Now assume that $G$ does not act transitively
on the complement of the diagonal in $\Lambda\times\Lambda$.
Let $\Gamma$ be an irreducible lattice in $H$
and let $\Gamma_0$ be the preimage of $\Gamma$
under the projection $G\to H=G/N$. We may assume that
$\Gamma_0$ contains a rank-one element. Moreover,
$\Gamma_0$ is non-elementary since this is the case
for $G$. Let $\mu$
be a (bi-invariant) Haar measure on the compact
normal subgroup $N$ and let
$\nu_0,\nu$ be a Haar measure on $\Gamma_0,\Gamma$. 
Then there is a continuous linear map
\[L^2(\Gamma_0,\nu_0)\to L^2(\Gamma,\nu)\] obtained
by mapping a square integrable 
function $f_0:\Gamma_0\to \mathbb{R}$ 
to the function $f:\Gamma\to \mathbb{R}$ defined by
\[f(g)=\int_Nf_0(gn)d\mu(n).\]
Note that this map is equivariant under the natural
left action of $G$. Via this map, the space
$L^2(\Gamma,\nu)$ is a coefficient module for $\Gamma_0$.

The subgroup $\Gamma_0$ of 
$G$ is closed and hence 
$H_{b}^2(\Gamma_0,L^2(\Gamma_0,\nu_0))\not=\{0\}$
by Proposition \ref{boundedcoho2}.
Composition with the coefficient map
$L^2(\Gamma_0,\nu_0)\to L^2(\Gamma,\nu)$ shows
that we have $H_b^2(\Gamma_0,L^2(\Gamma,\nu))\not=\{0\}$ as well
(see the discussion in \cite{M} for details).
This means that there are (not necessarily continuous)
unbounded twisted quasi-morphisms on $\Gamma_0$ 
with values in $L^2(\Gamma,\nu)$ which do not admit a
cocycle at bounded distance. By the remark following
Proposition \ref{boundedcoho2},
we may assume that the restriction of such a quasi-morphism
$\phi_0$ 
to each compact subset of $\Gamma_0$ is uniformly bounded.
By the defining inequality for a twisted quasi-morphism,
this implies that there is a universal constant
$c>0$ such that $\Vert \phi_0(h_1)-\phi_0(h_2)\Vert_p\leq c$
whenever $h_1,h_2$ project to the same element of $\Gamma$. 
Then the
quasi-morphism $\phi_0$ determines (non-uniquely) an unbounded
$L^2(\Gamma,\nu)$-valued twisted quasi-morphism 
$\phi$ on $\Gamma$ by defining $\phi(g)=\phi_0(h)$ for some
$h\in \Gamma_0$ which projects to $g$. This twisted quasi-morphism
is not at bounded distance from a cocycle since this
was not the case for $\phi_0$. Therefore 
we have 
\[H_{cb}^2(\Gamma,L^2(\Gamma,\mu))\not=0.\]
As before, the Burger-Monod super-rigidity
result for cohomology \cite{BM02} shows that
$H$ is simple of rank one (compare Theorem 14.2.2 of
\cite{M}).
\end{proof}

Now we are ready for the proof of the corollary from
the introduction (which is immediate from Corollary 1.24 of \cite{CM08}).
For this recall that a simply connected complete Riemannian
manifold $\tilde M$ of non-positive sectional curvature
is called \emph{irreducible} if $\tilde M$ does not split
as a non-trivial product. We have.

\begin{corollary}\label{rankrig}
Let $M$ be a closed Riemannian manifold of non-positive
sectional curvature. If the universal 
covering $\tilde M$ of $M$ is irreducible 
and if ${\rm Iso}(\tilde M)$ contains a parabolic element
then $M$ is locally symmetric. 
\end{corollary}
\begin{proof}
Let $M$ be a closed Riemannian manifold of non-positive 
sectional curvature with irreducible universal covering $\tilde M$.
The fundamental group $\pi_1(M)$ of $M$ acts cocompactly on 
the Hadamard space $\tilde M$ as a group of 
isometries. By the celebrated rank-rigidity theorem
(we refer to \cite{B95} for a discussion and for
references), either $\pi_1(M)$ contains a rank-one
element or $M$ is locally symmetric of higher rank.

Now assume that $\pi_1(M)$ contains a rank-one element.
By Lemma \ref{amenrad}, the amenable radical $N$ of ${\rm Iso}(\tilde M)$
fixes a point $x\in X$, and it fixes the limit
set of $\pi_1(M)$ pointwise. Since the 
action of $\pi_1(M)<{\rm Iso}(\tilde M)$ 
on $\tilde M$ is cocompact, the limit
set $\Lambda$ of $\pi_1(M)$ is the entire ideal
boundary $\partial \tilde M$ of $\tilde M$. Then
$N$ fixes
every geodesic ray issuing from $x$. 
This implies that $N$ is trivial.

By Theorem 1, either the isometry group of 
$\tilde M$ is an almost connected simple Lie group $G$ of rank one 
or ${\rm Iso}(\tilde M)$ is totally disconnected.
However, in the first case $\pi_1(M)$ is necessarily
a cocompact lattice in $G={\rm Iso}(\tilde M)$ since the action of 
${\rm Iso}(\tilde M)$ on $\tilde M$ is proper and 
cocompact.  Moreover, the dimension of the symmetric
space $G/K$ associated to $G$ coincides with the 
cohomological dimension of any of its uniform lattices and
hence it coincides with the dimension of $M$. 
But then the action of $G$ on $\tilde M$ is open.
Since this action is also closed, the action is transitive
and hence $\tilde M$ is a symmetric space.

We are left with showing that if ${\rm Iso}(\tilde M)$ contains
a parabolic element then the isometry group of 
$\tilde M$ is not totally disconnected.
Assume to the contrary that ${\rm Iso}(\tilde M)$ is totally 
disconnected. 
Since the action of ${\rm Iso}(\tilde M)$ on $\tilde M$ is 
cocompact, there is no non-trivial closed convex 
${\rm Iso}(\tilde M)$-invariant subset 
of $\tilde M$. Since ${\rm Iso}(\tilde M)$ is 
totally disconnected, by Theorem 5.1 of \cite{CM08},
point stabilizers of ${\rm Iso}(\tilde M)$ are open.
Then Corollary 3.3 of \cite{C07} implies that
every element of ${\rm Iso}(\tilde M)$ with vanishing
translation length is elliptic. This is a contradiction
to the assumption that ${\rm Iso}(\tilde M)$
contains a parabolic isometry.
\end{proof}

Finally we show Corollary 2 from the introduction.

\begin{corollary}
Let $G$ be a semi-simple Lie group with finite
center, no compact factors and rank at least 2. Let 
$\Gamma<G$ be an irreducible lattice, let $X$ be a proper
${\rm Cat}(0)$-space and let $\rho:\Gamma\to {\rm Iso}(X)$
be a homomorphism. If $\rho(\Gamma)$ is non-elementary and 
contains a
rank-one element then there is closed subgroup $H$  
of ${\rm Iso}(X)$ which is a compact extension of a simple Lie group
$L$ of rank one and there is a surjective
homomorphism $\rho:G\to L$.
\end{corollary}
\begin{proof}
Let $\Gamma<G$ be an irreducible lattice and let
$\rho:\Gamma\to {\rm Iso}(X)$ be a homomorphism.
Let $H<{\rm Iso}(X)$ be the closure of $\rho(\Gamma)$. 
If $\rho(\Gamma)$ is non-elementary and contains
a rank-one element then the same is true for $H$.
Then there is a non-trivial bounded cohomology class
in $H^2_{b}(H,L^2(H,\mu))$ which induces via $\rho$
a non-trivial bounded cohomology
class in $H^2_{b}(\Gamma,L^2(H,\mu))$.
Since $\Gamma$ is discrete, we have
$H_{cb}^2(\Gamma,L^2(H,\mu))\not=0$. 
By the constructions of Burger and Monod (see \cite{M}),
via inducing we deduce that the second bounded cohomology
group
$H_{cb}^2(G,L^{[2]}(G/\Gamma,L^2(H,\mu)))$ does not vanish,
where $L^{[2]}(G/\Gamma,L^2(H,\mu))$ denotes the
Hilbert $G$-module of all measurable maps $G/\Gamma\to 
L^2(H,\mu)$ with the additional property
that for each such map $\phi$ the function 
$x\to \Vert \phi(x)\Vert$ is square integrable on $G/\Gamma$
with respect to the projection of the Haar measure.
Here the $G$-action is determined by the homomorphism $\rho$.

By the results of Monod and Shalom \cite{MS06}, 
if $G$ is simple then there is a $\rho$-equivariant 
map $G/\Gamma\to L^2(H,\mu)$. 
However, this implies that $H$ is compact which is impossible
since $H$ contains a rank-one element.

If $G=G_1\times G_2$ for semi-simple Lie groups $G_1,G_2$
with finite center and no compact factor then the results
of Burger and Monod \cite{BM99,BM02} show that via possibly
exchanging $G_1$ and $G_2$ we may assume that there is 
a $G_1$-equivariant map $G/\Gamma\to L^2(H,\mu)$ for the
restriction of $\rho$ to $G_1$. Since $\Gamma$ is 
irreducible by assumption, the action of 
$G_1$ on $G/\Gamma$ is ergodic \cite{Z}. 
By Lemma \ref{amenrad}, the amenable radical $N$ 
of $H$ is compact and we 
deduce as in \cite{MS04} that there is a continuous homomorphism
$\psi:G\to H$. Since $G$ is connected, the image $\psi(G)=H/N$ is 
connected and hence by Proposition \ref{firstpart}, 
$H/L $ is a simple Lie group of rank one.
\end{proof}

\bigskip

\noindent
MATHEMATISCHES INSTITUT DER UNIVERSIT\"AT BONN\\
BERINGSTRA\SS{}E 1\\
D-53115 BONN\\

\smallskip

\noindent

\end{document}